\documentclass[a4paper,11pt]{amsart}

\usepackage{amsmath,latexsym,amsbsy,amssymb}
\usepackage{enumerate}
\usepackage{amsthm}
\usepackage[latin1]{inputenc}

\newtheorem{Theorem}{Theorem}[section]
\newtheorem{Definition}[Theorem]{Definition} 
\newtheorem{Proposition}[Theorem]{Proposition}
\newtheorem{Corollary}[Theorem]{Corollary}
\newtheorem{Rem}[Theorem]{Remark}
\newtheorem{Lemma}[Theorem]{Lemma}

\numberwithin{equation}{section}

\newcommand{\norm}[1]{\left\|#1\right\|}
\newcommand{\pd}[1]{\langle #1 \rangle}

\newcommand{\R}{{\mathbb R}}
\begin{document}

\title{On the minimum time function around the origin}

\author[Giovanni Colombo]{Giovanni Colombo}
\address[Giovanni Colombo]{Universit\`a di Padova, Dipartimento di Matematica Pura ed Applicata, via Trieste 63, 35121 Padova, Italy}
\email{colombo@math.unipd.it}

\author[Khai T. Nguyen]{Khai T. Nguyen}
\address[Khai T. Nguyen]{Universit\`a di Padova, Dipartimento di Matematica Pura ed Applicata, via Trieste 63, 35121 Padova, Italy}
\email{khai@math.unipd.it}

\thanks{This work was partially supported by M.I.U.R., project
``Viscosity, metric, and control theoretic methods for nonlinear
partial differential equations''  and by CARIPARO}



\date{\today}
\begin{abstract}
We deal with finite dimensional linear and nonlinear control systems. If the system is linear and autonomous and satisfies the 
classical \textit{normality} assumption, we improve the well known result on the strict convexity of the reachable set from the origin by giving a
polynomial estimate. The result is based on a careful analysis of the switching function.
We extend this result to nonautonomous linear systems, provided the time dependent system is not too far from
the autonomous system obtained by taking the time to be $0$ in the dynamics.\\
Using a linearization approach, we prove a bang-bang principle, valid in dimensions $2$ and $3$ for a class of nonlinear systems, affine and symmetric
with respect to the control. Moreover we show that, for two dimensional systems, the reachable set from the origin
satisfies the same polynomial strict convexity property as for the linearized dynamics, provided the nonlinearity is small enough.
Finally, under the same assumptions we show that the epigraph of the minimum
time function has positive reach, hence proving the first result of this type in a nonlinear setting. 
In all the above results, we require that the linearization at the origin be normal.\\
We provide examples showing the sharpness of our assumptions.
\end{abstract}
\maketitle
\section{Introduction}
In the theory of autonomous linear control systems, the assumption of \textit{normality}, i.e., a strong controllability assumption requiring that
if each control component is used separately, then Kalman rank condition is satisfied, is well known. In particular (see \cite[Sections 14, 15, and 16]{hls}),
normality implies that the control steering the origin to a point $x$ in minimum time is unique and bang-bang; moreover the reachable set
from the origin, $\mathcal{R}^\tau$, is a \textit{strictly} convex body for all times $\tau >0$. Simple examples, on the other hand,
show that convexity of the reachable set is easily lost when passing
to a nonlinear dynamics, even if the control covers all directions and appears linearly. In particular, the control system
\begin{equation}\label{eq:example_intro}
\begin{cases}
\dot{y}_1(t)&= - [y_2(t)]^2 + u_1 (t),\quad -1\le u_1 \le 1\\
\dot{y}_2(t)&= u_2(t) \qquad\qquad\qquad\, -1\le u_2 \le 1
\end{cases}                                            
\end{equation}
which was analyzed in \cite{CS0}, fails to have convex (and even normally regular, i.e., ``without inward corners'')
reachable sets from the origin at any positive time.\\
Similarly, the minimum time to reach the origin for a linear system is semiconvex, i.e., it is a quadratic perturbation of a convex function,
provided a first order controllability assumption is satisfied (see \cite{CS0,CS}), or has epigraph with positive reach, i.e., it satisfies a quite
good kind of generalized convexity, provided Kalman rank condition holds (see \cite{CM,CMW}). On the other hand, the same example introduced above shows
that such properties may fail even for a simple nonlinear dynamics.\\
We are not aware of any result of semiconvexity type valid for a nonlinear dynamics, where actually results of semiconcave
type are more natural and easier
to obtain, provided the target (or the dynamics) satisfy an inner ball condition (see \cite{CS,cafra,CK,K,camapw,CaKh}).
Further results on the regularity
of the minimum time function $T$, for two dimensional systems with a single input control, are described in \cite[Chapter 3]{BP},
where in particular, under generic assumptions, a characterization of smooth and nonsmooth points of the level sets of $T$ is given.\\
In this paper, we give a contribution to the understanding of the behavior of the minimum time function both for linear,
autonomous and nonautonomous, and nonlinear control systems. We take the origin as the target (or as the source, for the reversed dynamics),
hence going outside the realm of semiconcavity, and give results, under easily verifiable assumptions, on the following topics:
\begin{itemize}
\item strict convexity of the reachable set from the origin
\item uniqueness of the optimal control
\item a nonlinear bang-bang principle
\item extending backwards optimal trajectories
\item positive reach of the epigraph of $T$.
\end{itemize}
Our method is based on new results on linear control systems satisfying the normality condition, and on linearization at the origin.
The underlying idea, in fact, is requiring enough strength to such linearization,
so that examples like \eqref{eq:example_intro} are ruled out.
Our linear results hold in any space dimension, while the nonlinear part is confined to two or three dimensional
spaces. Our arguments are based essentially on a careful analysis of the \textit{switching function}, namely the function whose
sign is expected to determine the optimal control, according to Pontryagin's Maximum Principle.
The point is exactly showing that this sign is well defined,
except for at most finitely many zeros. To this aim, the normality assumption is pivotal, as it permits to split any finite interval into finitely many sets,
each one being a disjoint union of finitely many intervals, where the switching function or its derivatives are uniformly bounded away from zero
(see Lemma \ref{7fin int}). From this fact we are able to deduce a quantitative estimate on the strict convexity of the reachable set. More precisely
(see Theorem \ref{7LA}), for a linear control system in $\mathbb{R}^N$ we show that for all $\tau >0$ there exists a positive constant $\gamma >0$ such that 
\[
 \langle \zeta , y-x \rangle \le - \gamma \| \zeta\|\,\| y-x \|^N\qquad \text{for all } x,y\in \mathcal{R}^\tau,\;\zeta\in N_{\mathcal{R}^\tau}(x)
\]
(here $N_{\mathcal{R}^\tau}(x)$ denotes the normal cone to $\mathcal{R}^\tau$ at $x$). We show through an example that the exponent $N$ is optimal.
Section \ref{sec:lin} is devoted to the above topic, together with an auxiliary study for a linear nonautonomous dynamics.\\
The nonlinear part starts with a nonlinear bang-bang result (Section \ref{sec:bangbang}), valid up to dimension $3$.
We consider a nonlinear control system which is affine with respect to the control: if the linearization at the origin is normal, then every optimal
control is bang-bang. The proof is based on Pontryagin's Maximum Principle: if the nonlinearity contains only parts which are of order larger or equal
to the space dimension, then we are able to transfer to the switching function all the properties satisfied by the switching function of the dynamics
linearized at the origin. This idea is at the basis also of the strict convexity of the reachable set for a nonlinear two dimensional dynamics
(see Theorem \ref{7NL}) and of proving that all points close enough to the origin are \textit{optimal}, i.e., any trajectory steering a point to the origin
optimally can be extended backwards still remaining optimal (see Theorem \ref{7optimal}). In this case, the difficulty is extending the optimal control: our analysis
permits to predict backwards the sign of the switching function. 
Finally, we show under the same assumptions that the epigraph of the minimum time function $T$ has positive reach, hence obtaining
a rich bunch of regularity properties for $T$, as listed in Theorem \ref{diff}.
We also show through an example (see Example \ref{exnonstrict}) that the assumptions on the nonlinear
part cannot be avoided, while \eqref{eq:example_intro} shows that normality at the origin is essential.
The restrictions on the space dimension for the nonlinear results will be explained after the relevant proofs (see Remarks \ref{dimrestr1} and \ref{dimrestr2}).\\
To our knowledge, the results we present here as well as most of the used methods do not trace back to previous literature.
In particular, the nonlinear bang-bang results of Krener \cite{kre} and Sussmann \cite{suss} seem to be of a very different nature.
\section{Preliminaries}\label{sec:prelim}
\subsection{Nonsmooth analysis and sets with positive reach}\quad\\
In all the paper, the space dimension $N$ will be supposed larger or equal to $2$.\\
Let $K\subset\R^N$ be closed. The distance function to $K$ and the projection mapping onto $K$ are defined respectively by
\begin{eqnarray*}
d_K(x)&=&\min\{\|y-x\|\, :\,y\in K\}, \\
\pi_K(x)&=&\{y\in K\, :\,\|y-x\|=d_K(x)\}.
\end{eqnarray*}
The boundary of $K$ is $\mathrm{bdry}\, K$. Given $x\in K$ and
$v\in\mathbb R^N$, we say that $v$ is a \emph{proximal normal} to $K$ at $x$ (and we will denote this fact by $v\in N_K(x)$)
if there exists $\sigma=\sigma(v,x)\ge 0$ such that
\[
\pd{  v,y-x}\le\sigma \norm{y-x}^2\quad\textrm{ for all }y\in K;
\]
equivalently  
$v\in N_K(x)$ if and only if
there exists $\rho >0$ such that $\pi_K(x+\rho v)=\{x\}$, and in this case we say that the proximal normal $v$ 
\textit{is realized by a ball of radius $\rho>0$}. 
If $K$ is convex then $N_K(x)$ coincides with the normal cone of Convex Analysis.
The proximal subdifferential of a function $f:\Omega\rightarrow \R$ at a point $x$ of its domain, $\partial f(x)$, is the set of vectors $v\in \R^N$
such that 
\[
(v,-1)\in N_{\text{epi}(f)}(x,f(x)),
\]
where epi$(f)$, the epigraph of $f$, is defined as
\[
\text{epi}(f)=\{ (x,y)\in\Omega\times\R : y\ge f(x)  \}.
\]
For an introduction to nonsmooth analysis we make reference, e.g., to \cite[Chapters 1 and 2]{CLSW}.\\
We will make use of the following concepts due to Federer \cite{federer}: given an arbitrary set $K\subset \R^N$, we set
\begin{eqnarray*}
\text{Unp}(K) &=& \{ x\in\mathbb R^N\, :\, \pi_{K}(x)\text{ is a singleton}\}, \\
\text{reach}(K,x) &=&\sup \{ r\ge 0 \, :\, x+r\mathbb B^n\subseteq\text{Unp}(K)\},
\end{eqnarray*}
the latter being defined for $x\in K$. 
We remark that $\text{reach}(K,x)$ is continuous with respect to $x\in K$.
\begin{Definition}
We say that a closed $K\subseteq\mathbb R^N$ has \emph{positive reach} if 
$$\mathrm{reach}(K):=\inf_{x\in K}\mathrm{reach}(K,x)>0.$$
A locally closed set $K$ has \emph{locally positive reach} if $\mathrm{reach}(K,x)>0$ for every $x\in K$.
\end{Definition}
The positive reach property is actually an external sphere condition with (locally) uniform radius. More precisely, it holds
(see \cite[\S 4]{federer}):
\begin{Proposition}\label{distproj} 
Let $K\subset \R^N$ be closed. Then $K$ has positive reach if and only if there
exists a continuous function $\varphi:K\to[0,+\infty)$ such that the inequality
\begin{equation}\label{ineqset}
\langle v,y-x\rangle\le\varphi(x)\| v\|\,\| y-x\|^2
\end{equation}
holds for all $x,y\in K$ and $v\in N_{K}(x)$.
\end{Proposition}
We say that $f:\Omega\rightarrow \mathbb{R}$ is of class $\mathcal{C}^{1,1}(\Omega)$ (here $\Omega\subseteq\mathbb{R}^N$ is open) if its partial derivatives exist
and are Lipschitz.
\begin{Rem}
The case when $\varphi(x)\equiv 0$ in \eqref{ineqset} is equivalent to the convexity of $K$, and \eqref{ineqset} is in this sense a generalization of convexity.
Furthermore, it is easy to see that if the boundary of $K$ is the graph of a $\mathcal{C}^{1,1}$-function, then
$K$ has positive reach and $\varphi$ in \eqref{ineqset} can be taken as constant. Hence positive reach property generalizes $\mathcal{C}^{1,1}$-manifolds as well.
\end{Rem}
Lower semicontinuous functions whose epigraph has positive reach enjoy remarkable regularity properties, which
are similar to properties of convex functions. In particular, the following result holds true (we state it for continuous functions
for simplicity). The Lebesgue $N$-dimensional measure and the Hausdorff $d$-dimensional measure are denoted,
respectively, by $\mathcal{L}^N$ and $\mathcal{H}^d$.
\begin{Theorem}\label{diff}
Let $\Omega\subset\mathbb R^N$ be open, and let
$f:\Omega\to\mathbb R$ be continuous, and such that
$\text{epi}(f)$ has locally positive reach.
Then there exists a sequence of sets $\Omega_h\subseteq\Omega$ such
that $\Omega_h$ is compact in $\Omega$ and
\begin{enumerate}
\item[(1)] the union of $\Omega_h$ covers $\mathcal L^N$-almost all $\Omega$;
\item[(2)] for all $x\in\bigcup_h\Omega_h$ there exist $\delta=\delta(x)>0$, $L=L(x)>0$ such that
\[
f\textrm{ is Lipschitz on $B(x,\delta)$
with ratio }L, \textrm{and hence semiconvex on }B(x,\delta).
\]
\end{enumerate}
Consequently,
\begin{enumerate}
\item[(3)] $f$  is a.e. Fr\'echet differentiable and admits a second order Taylor expansion
around a.e. point of its domain.
\end{enumerate}
Moreover, the set of points where the graph of $f$ is nonsmooth has small Hausdorff dimension. More precisely, 
\begin{enumerate}
\item[(4)] for every
$k=1,\ldots ,N$, the set $\{ x\in \mathrm{int\, dom}(f)\ |\ \text{the dimension of }$
$ \partial f(x)\text{ is } \ge k\}$
is countably $\mathcal{H}^{N-k}$-rectifiable.
\end{enumerate}
Finally,
\begin{enumerate}
\item[(5)] $f$ has locally bounded variation in $\Omega$.
\end{enumerate}
\end{Theorem}
This result is essentially Theorem 5.1 and Proposition 7.1 in \cite{CM}. For properties of semiconvex/semiconcave functions we
refer to \cite{CS}.
\subsection{Control theory}
We will consider control systems linear or nonlinear with respect to the space variable and affine and symmetric
with respect to the control. More precisely, we will consider the linear control system 
\begin{equation}\label{7LNsystem}
\left\{\begin{array}{ll}
\dot{y}(t)\: = \: Ay(t)+Bu(t) & \text{a.e.}\\
u(t) \: \in \:  \mathcal{U} =[-1,1]^{M}& \text{a.e.}\\
y(0)  \: = \:  0,
\end{array}\right.
\end{equation}
where $1\le M\le N$ and $A\in\mathbb{M}_{N\times N}$, $B\in\mathbb{M}_{N\times M}$,
being possibly time dependent, and $\mathcal{U}=[-1,1]^{M}\ni (u_1,\ldots,u_M)=:u$,
together with the nonlinear control system
\begin{equation}\label{7NLNsystem}
\left\{\begin{array}{ll}
\dot{y}(t)\: = \: F(y(t))+G(y(t))u(t) & \text{a.e.}\\
u(t) \: \in \:  [-1,1]^{M} & \text{a.e.}\\
y(0)  \: = \:  0,
\end{array}\right.
\end{equation}
where $F$ and $G$ are suitable vector fields (the actual assumptions will be stated later).
We will use also the notation $B=(b_1,\ldots ,b_M)$ or $G=(g_1,\ldots,g_M)$, where each entry is an $N$-dimensional column.
We denote by $\mathcal{U}_{ad}$, the set of admissible controls, i.e., all measurable
functions $u$, such that $u(s)\in\mathcal{U}$ for a.e. $s$.
For any $u(\cdot)\in\mathcal{U}_{ad}$, the (unique, as it will follow from the assumptions on $F$ and $G$)
Carath\'eodory solution of (\ref{7LNsystem}) or of \eqref{7NLNsystem}
is denoted by $y^{u}(\cdot)$. In the linear case,
\[
y^{u}(t)=\int_0^te^{A(t-s)}Bu(s)ds,
\]
so that the reachable set from $0$ in time $t$ can be described by 
\begin{equation}\label{reprRT}
\mathcal{R}^t=\Big\lbrace{\int_0^te^{A(t-s)}Bu(s)ds)\ |\ u(\cdot)\in\mathcal{U}_{ad}\Big\rbrace}.
\end{equation}
It is well known that in the linear case the set $\mathcal{R}^t$ is convex and compact (see, e.g., \cite[Lemma 12.1]{hls}),
while in the  nonlinear case \eqref{7NLNsystem} $\mathcal{R}^t:=\lbrace{ y^u(t)\ |\ u(\cdot)\in\mathcal{U}_{ad}\rbrace}$
is compact and not necessarily convex (see, e.g., \cite[Chapter 10]{cesari}).\par
For a fixed $x\in\mathbb{R}^N$, we define
$$\theta(x,u):=\min\ \lbrace{t\geq 0\: |\: y^{x,u}(t)=0\rbrace},$$
where $y^{x,u}(\cdot)$ denotes the solution of $\dot{y}=F(y)+G(y)u$ such that $y(0)=x$.
Of course, $\theta(x,u)\in [0,+\infty]$, and $\theta(x,u)$ is the time taken for the trajectory
$y^{x,u}(\cdot)$ to reach $0$, provided $\theta(x,u)<+\infty$.
The \textit{minimum time} $T(x)$ to reach $0$ from $x$ is defined by
\begin{equation}\label{minT}
T(x):=\inf\ \lbrace{\theta(x,u)\: |\: u(\cdot)\in\mathcal{U}_{\mathrm{ad}}\rbrace}.
\end{equation}
Observe that the sublevel $\mathcal{R}_t=\{ x: T(x)\le t\}$ of $T(\cdot)$ equals the reachable
set from the origin within the (same) time $t$ for the \textit{reversed dynamics} 
\[
\dot{x}=-F(x)-G(x)u,\quad u\in\mathcal{U}. 
\]
If $\bar{u}$ is an admissible control steering $x$ to the origin in the minimum time $T(x)$ (i.e., an \textit{optimal}
control), then the Dynamic Programming Principle (see, e.g., Proposition 2.1, Chapter IV, in \cite{BCD})
implies that $T(\cdot)$ is strictly increasing along the optimal trajectory $y^{x,\bar{u}}$. Therefore, for all
$0<t<T(x)$ the point $y^{x,\bar{u}}(t)$ belongs to the boundary of $\mathcal{R}_{t}$.

Pontryagin's Maximum Principle is a fundamental tool for the analysis of optimal control problems. 
We state it for points belonging to the boundary of reachable sets. In view of the previous remark, this will apply also to
points belonging to optimal trajectories. We give first its linear version.
\begin{Theorem}[Pontryagin's Maximum Principle for linear systems] \label{th:PMPlin}
Consider the problem \eqref{7LNsystem}, fix $T>0$, and suppose $\bar{x}\in\mathcal{R}^T$
is realized by the control $\bar{u}(\cdot)\in\mathcal{U}_{ad}$ (i.e., $y^{\bar{u}}(T)=\bar{x}$).
Then $\bar{x}\in\mathrm{bdry}\,\mathcal{R}^T$ if and only if for some
$\zeta\in N_{\mathcal{R}^T}(\bar{x})$, $\zeta\neq 0$, it holds
\begin{equation}\label{signlin}
\bar{u}_i(t)=\mathrm{sign} \langle \zeta,e^{A(T-t)}b_i\rangle,\quad a.e.\ t\in [0,T],
\end{equation}
for all $i=1,2,\ldots ,M$.
\end{Theorem}
A well known reference for this result is \cite[Lemma 13.1]{hls}.

Before stating Pontryagin's principle for the nonlinear case \eqref{7NLNsystem} (with $M=1$),
we need to introduce the \textit{Maximized Hamiltonian}.\\
We define for every triple $(x,p,u)\in\R^N\times\R^N\times [-1,1]$
\[
 \mathcal{H}(x,p,u)=\langle p, F(x)\rangle + u\langle p,G(x)\rangle
\]
and
\[
H(x,p)=\max\{ \mathcal{H}(x,p,u) : u\in [-1,1] \}. 
\]
Then Pontryagin's principle reads as follows (see, e.g., \cite[Section 2.1]{BP} or \cite[Theorem 3.5.4]{clarke}).
\begin{Theorem}[Pontryagin's Maximum Principle for nonlinear systems]\label{th:PMPnonlin}
Fix $T>0$ and let $y^u$ be a trajectory of \eqref{7NLNsystem} such that $y^u(T)$ belongs to the boundary of
$\mathcal{R}^T$. Then there exist a function $\lambda : [0,T]\rightarrow\R^N$, never vanishing, and a constant
$\lambda_0 \le 0$ such that for a.e. $t\in [0,T]$ one has:
\begin{itemize}
\item[i)] $\dot{\lambda}(t)=-\lambda(t) \big(DF(y^u(t))+DG(y^u(t))u(t)\big)$,
\item[ii)] $\mathcal{H}(y^u(t),\lambda(t),u(t))+\lambda_0=0$, 
\item[iii)] $\mathcal{H}(y^u(t),\lambda(t),u(t))=H(y^u(t),\lambda(t))$.
\end{itemize}
Furthermore,
\[
\lambda(T)\in N_{\mathcal{R}^T}(y^u(T))
\]
where here $N_{\mathcal{R}^T}(y^u(T))$ denotes the \textit{Clarke} normal cone.
\end{Theorem}
\begin{Definition}\label{def:essdet}
We say that a control $u$ is essentially determined by Pontryagin's Principle if for any $u_1$ satisfying iii) in Theorem
\ref{th:PMPnonlin} (for the adjoint curve $\lambda$ associated with the trajectory $y^u$)
one has $u_1(t)=u(t)$ a.e. in $[0,T]$.
\end{Definition}
In the following we will make extensive use of the classical concept of normality for linear systems, which we are now going to introduce.
\begin{Definition}
The system (\ref{7LNsystem}) is normal if and only if for every column $b_i$ of $B$, $i=1,\ldots ,M$, we have
\[
\mathrm{Rank}\big[b_i,Ab_i,\ldots ,A^{N-1}b_i\big]=N.
\]
\end{Definition}
The main classical result for normal linear systems is concerned with the reachable set.
\begin{Theorem}\label{7Strict-weak}
Assume that the linear control system(\ref{7LNsystem}) is normal. Then the reachable set $\mathcal{R}^T$ is strictly convex for any $T>0$.
\end{Theorem}
\textbf{Proof.} One can find a proof in \cite{hls}, Sections 14 and 15.
\qed
\begin{Rem}
If the system (\ref{7LNsystem}) is normal then $(A,B)$ satisfies the Kalman rank condition. Therefore the minimum time function
is everywhere finite and continuous (actually H\"older continuous with exponent $1/N$, see, e.g., Theorem 17.3 in \cite{hls} and Theorem 1.9, Chapter IV,
in \cite{BCD} and references therein).
\end{Rem}
\section{Quantitative strict convexity of reachable sets: the linear case}\label{sec:lin}
\subsection{Autonomous systems}\quad\\
This subsection is devoted to improving the classical result on the strict
convexity of reachable sets for normal linear control systems of the type \eqref{7LNsystem}. We will give an estimate for the boundary of reachable
sets which implies a uniform (polynomial) strict convexity with an optimal exponent.\\
In the first Lemma we define the switching function and begin studying its behavior under the normality assumption.
\begin{Lemma}\label{7kahn}
Let $A\in\mathbb{M}_{N\times N}$ and $b\in\mathbb{R}^N$ be such that 
\begin{equation}\label{7rank cond}
\mathrm{Rank}\big[b,Ab,\ldots ,A^{N-1}b\big]=N.
\end{equation}
Take $\zeta\in\mathbb{R}^N$, with $\|\zeta\|=1$, and define, for $s\in [0,+\infty)$
\begin{equation}\label{7defg}
g(s)\: =\: \langle e^{As}b\ ,\ \zeta\rangle.
\end{equation}
Then there exists a constant $\mathcal{L}$, depending only on $A,b,N$ such that, for all $s\in [0,+\infty)$,
\begin{equation}\label{7sumder}
\sum_{i=0}^{N-1}|g^{(i)}(s)|\geq \mathcal{L}e^{-\|A\|s}.
\end{equation}
\end{Lemma}
\textbf{Proof.}
Set 
\begin{eqnarray*}
K=\begin{pmatrix} 
b,\;\,
Ab,\;\,
\ldots,\;\,
A^{N-1}b
\end{pmatrix}
\end{eqnarray*}
and observe that, by (\ref{7rank cond})
\begin{equation}\label{7lbound}
\mathcal{L}=\min_{\|\zeta\|=1}\ \| K\zeta\|\ >0.
\end{equation}
Fix $\zeta\in\mathbb{R}^{N}$ with $\|\zeta\|=1$ and write $\zeta_1(s)=e^{sA^T}\zeta$. Observe
that $\zeta=e^{-A^{T}s}\zeta_1(s)$ and $\|\zeta_1(s)\|\geq e^{-s\|A\|}$. We compute now, for $i=0,1,\ldots ,N-1,$
\begin{equation}\label{compgi}
g^{(i)}(s)=\langle e^{As}A^ib,\zeta\rangle=\langle A^ib,\zeta_1(s)\rangle.
\end{equation}
Therefore,
\begin{eqnarray*}
K\zeta_1(s)&=&\begin{pmatrix} 
b,\;\,
Ab,\;\,
\ldots,\;\,
A^{N-1}b
\end{pmatrix}
\zeta_1(s)\\
&=&\begin{pmatrix}
g^{(0)}(s),\;\,
g^{(1)}(s),\;\,
\ldots,\;\,
g^{(N-1)}(s)
\end{pmatrix}.
\end{eqnarray*}
Using (\ref{7lbound}) we have that 
\[
\|K\zeta_1(s)\|\geq\mathcal{L}e^{-s\|A\|}.
\]
On the other hand,
\[
\|K\zeta_1(s)\|\leq\sum_{i=0}^{N-1}|g^{(i)}(s)|
\]
and the proof is concluded. \qed

\medskip

The next Lemma is crucial for estimating the number of zeros of the switching function $g$
(corresponding to the number of switching points of the optimal control associated with $g$) and for studying their multiplicity.
We recall that the constant $\mathcal{L}$ was defined in \eqref{7lbound}.
\begin{Lemma}\label{7fin int} 
Let $A\in\mathbb{M}_{N\times N}$ and $b\in\mathbb{R}^N$ be satisfying (\ref{7rank cond}). Take $\zeta\in\mathbb{R}^N$, $\|\zeta\|=1$,
and fix $T>0$. Let $g(s)$, $s\in [0,T]$, be defined as in (\ref{7defg}).\\
\indent Then there exist disjoint sets $I_0,\ldots ,I_{N-1}$ and numbers $\mathcal{N}_i$, depending only on $A,b,T$ and $N$ such that 
\[
[0,T]=\bigcup_{i=0}^{N-1}I_i
\]
and, for all $i=0,1,\ldots ,N-1$, the set $I_i$ is the disjoint union of at most $\mathcal{N}_i$ intervals. Moreover, for each $i=0,1,\ldots ,N-1$,
for all $s\in I_i$, we have 
\begin{equation}\label{7estgi}
|g^{(i)}(s)|\geq\frac{\mathcal{L}}{N}e^{-\|A\|s}.
\end{equation}
\end{Lemma}
\textbf{Proof.}  We proceed by induction for $i$ from $0$ to $N-1$. Set 
\begin{equation}\label{defc}
c(s)=\frac{\mathcal{L}e^{-\|A\|s}}{N},
\end{equation}
and 
\[
J_0=\lbrace{s\in (0,T)\ |\ |g(s)|< c(s)\rbrace}.
\]
Since $J_0$ is open, we can write it as the disjoint union of at most countably many open intervals,
\begin{equation}\label{7J0}
J_0=\bigcup_{k=1}^{\infty}J_0^k.
\end{equation}
We assume, without loss of generality, that there are at least $N$ such intervals. Fix now any number $N' \ge N$,
and take a subfamily of the intervals
$J_0^k$ consisting of at most $N'$ elements. Without loss of generality, we can rearrange their indexes $k$ so that
$J_0^k=(a_{2k},a_{2k+1})$, where $1\le k\le N'$ and $0\le a_2<a_3\le a_4<a_5\le \ldots \le a_{2k}<a_{2k+1}\ldots <a_{2N'+1}\le T$.\\
\indent Now, fix $k$ and consider the $N$ intervals $(a_{2k},a_{2k+1}),\ldots ,(a_{2(k+N-1)},a_{2(k+N)-1})$. Set, for $j=0,1,\ldots ,$ $N-1$,
\[
(a_{2(k+j)},a_{2(k+j)+1}):=I^-_j,
\] 
and, for $j=0,1,\ldots ,N-2$
\[
[a_{2(k+j)+1},a_{2(k+j+1)}]:=I^+_j.
\]
Observe that 
\[
\bigcup_{j=0}^{N-1}I_j^-\cup\bigcup_{j=0}^{N-2}I_j^+=(a_{2k},a_{2(k+N)-1}).
\]
We are going to give a lower bound on $|a_{2(k+N)-1}-a_{2k}|$ which will turn out to be independent of both $k$ and $N'$.
From this fact it will follow automatically that the intervals
$(a_{2k},a_{2k+1})$ are nonempty only for finitely many $k$.\\
\indent Observe that for each $j=0,1,\ldots ,N-2$, there exists at least one point $c_j^1\in I_j^+$ such that $g'(c_j^1)=0$.
Therefore, there exist at least $N-2$ points, say $c^2_j$, for $j=0,1,\ldots ,N-3$, such that 
\[
c_j^2\in (c^1_{j},c^1_{j+1})\quad \mathrm{and}\quad g''(c_j^2)=0.
\]
Proceeding by induction we see that, for each $i=1,\ldots ,N-1$, there exists at least one point $c_i\in (a_{2k+1},a_{2(k+N)-1})$
such that $g^{(i)}(c_i)=0$. \\
Pick any $s_0\in (a_{2k},a_{2k+1})$. We have 
\begin{equation}\label{7g0}
|g(s_0)|<c(s_0),
\end{equation}
and, for $i=1,\ldots , N-1$,
\begin{eqnarray*}
|g^{(i)}(s_0)|&=&|g^{(i)}(s_0)-g^{(i)}(c_i)|=\Big|\int_{s_0}^{c_i}g^{(i+1)}(s)ds\Big|\\
&\leq&\int_{a_{2k}}^{a_{2(k+N)-1}}|g^{(i+1)}(s)|ds\leq (a_{2(k+N)-1}-a_{2k})e^{\|A\| T}\| A^{i+1}b\|,
\end{eqnarray*}
where the last inequality is due to \eqref{compgi}. Therefore,
\begin{equation}\label{7sumgi}
\sum_{i=1}^{N-1}|g^{(i)}(s_0)|\leq (a_{2(k+N)-1}-a_{2k})e^{\|A\| T}\sum_{i=1}^{N-1}\|A^{i+1}b\|.
\end{equation}
On the other hand, recalling (\ref{7sumder}), (\ref{defc}), and (\ref{7g0}) we have 
\begin{eqnarray}\label{7gn-1}
\sum_{i=1}^{N-1}|g^{(i)}(s_0)|\geq  \mathcal{L}e^{-\|A\|s_0}-c(s_0)=\frac{N-1}{N}\mathcal{L}e^{-\|A\|s_0}\geq \frac{N-1}{N}\mathcal{L}e^{-\|A\|T}.
\end{eqnarray}
From (\ref{7sumgi}) and (\ref{7gn-1}) we obtain
\begin{equation}\label{7length}
a_{2(k+N)-1}-a_{2k}\geq\frac{(N-1)\mathcal{L}e^{-2\|A\|T}}{N\sum_{i=1}^{N-1}\|A^{i+1}b\|},
\end{equation}
which is the desired estimate. Observe that the right hand side of (\ref{7length}) depends only on $A,b,T$ and $N$.\\
\indent We set now $\mathcal{N}_0$ to be the number of nonempty intervals contributing to the union in (\ref{7J0}), and recall that we have
just proved that $\mathcal{N}_0$ depends only on $A,b,T$ and $N$, and actually
\begin{equation}\label{7boundN0}
\mathcal{N}_0\leq\frac{N^2}{N-1}\frac{T}{\mathcal{L}}e^{2\|A\| T}\sum_{i=1}^{N-1}\|A^{i+1}b\| +N-1.
\end{equation}
Set $I_0=[0,T]\backslash J_0$ and observe that we have completed the proof of the lemma for $i=0$.\\
\indent After this step, we formulate our induction process. We are going to construct, for each $i=0,\ldots ,N-1$, two
disjoint sets $I_i,J_i$ with the following properties
\begin{enumerate}
\item [(ind1)] for every $s\in I_i$, $|g^{(i)}(s)|\geq c(s)$;
\item [(ind2)] for every $s\in J_i$, $\sum_{j=i+1}^{N-1}|g^{(j)}(s)|\geq (N-i-1)c(s)$;
\item [(ind3)] $J_i\bigcup I_i=J_{i-1}$;
\item [(ind4)] $J_i$ is a finite union of open intervals whose number is at most $\mathcal{N}_i$, and $\mathcal{N}_i$
depends only on $T,\mathcal{L},A,b,N,i$;
\item [(ind5)] $I_i$ is the finite union of at most $\mathcal{N}_i+\mathcal{N}_{i-1}$ intervals.
\end{enumerate}
For $i=0$ the above construction was already performed (take $J_{-1}=(0,T)$). Pick any $i=1,\ldots ,N-2$
(the case $i=N-1$ will be treated separately) and assume that (ind1), \ldots, (ind5) hold up to $i-1$. We wish to show that
the above statements hold for $i$ as well. To this aim, consider the set 
\[
J_{i}:=\lbrace{s\in J_{i-1}\ |\ |g^{(i)}(s)|<c(s)\rbrace}.
\]
For every connected component $(a,b)$ of $J_{i-1}$, we are going to prove that $J^{(a,b)}_{i}:=J_{i}\cap (a,b)$ is a
finite union of intervals, and give a bound on their number $\mathcal{N}_{i}^{(a,b)}$.\\
So, fix a connected component $(a,b)$ of $J_i$ and represent the open set $J^{(a,b)}_{i}$ as a disjoint union of at most
countably many intervals $J^k_{i}$, $k\in\mathbb{N}$ (for simplicity of writing we drop the dependence on $(a,b)$).
Assume without loss of generality that there are at most $N-i$ intervals $J^k_{i}$, fix any number $N'' \ge N-i$, and
take any subfamily of $\{ J^k_{i}\}$ consisting of at most $N''$ intervals. We can write $J^k_{i}= ( a_{2k}, a_{2k+1})$, where
$1\le k \le N''$ and $a\le a_2<a_3\le a_4<a_5\le\ldots < a_{2N''+1}\le b$. Fix $k$ and consider the intervals $(a_{2k},a_{2k+1}),\ldots ,
(a_{2(k+N-i-1)},a_{2(k+N-i)-1})$. Set, for $j=0,\ldots ,N-i-1$,
\[
(a_{2(k+j)},a_{2(k+j)+1}):= I_j^-, 
\]
and, for $j=0,\ldots ,N-i-2$,
\[
[ a_{2(k+j)+1}, a_{2(k+j+1)}]:=I_j^+.
\]
Observe that
\[
\bigcup_{J=0}^{N-i-1}I_j^-\cup\bigcup_{j=1}^{N-i-2} I_j^+= (a_{2k},a_{2(k+N-i)-1}).
\]
For each $j=0,\ldots ,N-i-2$ there exists at least one point $c_j^i\in I_j^+$ such that $g^{(i)}(c^i_j)=0$. Proceeding by induction
we see that for each $m=0,\ldots ,N-i-2$ there exists a point $c_m\in (a_{2k+1},a_{2(k+N-i)-1})$ such that $g^{(i+m)}(c_m)=0$.
Pick any $s_0\in (a_{2k},a_{2(k+1)})$. By arguing as for $i=0$, we obtain on one hand
\begin{equation}\label{7gi}
|g^{(i)} (s_0)| < c(s_0) 
\end{equation}
and, for all $m=0,\ldots N-i-1$,
\[
|g^{(i+m)}(s_0)|\le (a_{2(k+N-i)-1}-a_{2k})e^{\| A\|T}\big\| A^{i+m+1}b\big\|, 
\]
the latter inequality being due to \eqref{compgi}.
Thus
\begin{equation}\label{7sumgii}
\sum_{m=i+1}^{N-1} | g^{(m)} (s_0)|\le  (a_{2(k+N-i)-1}-a_{2k})e^{\| A\|T}\sum_{m=i+1}^{N-1}\big\| A^{m+1}b\big\|.
\end{equation}
On the other hand, owing to (ind2) and \eqref{7gi}, we obtain
\[
\sum_{m=i+1}^{N-1}|g^{(m)}(s)|\geq (N-i-1)c(s_0).
\] 
By combining the above inequality with \eqref{7sumgii} we now obtain
\[
a_{2(k+N-i)-1}-a_{2k} \ge \frac{N-i-1}{N}\mathcal{L} e^{-\| A\|T} \frac{1}{\sum_{m=i+1}^{N-1}\| A^{m+1}b\|}.                                                            
\]
Therefore, $J_{i}^{(a,b)}$ is the union of finitely many disjoint open intervals
$\big(a_{2k}^{i+1},a_{2k+1}^{i+1}\big)$, $k=1,\ldots ,\mathcal{N}^{(a,b)}_{i}$, where
\begin{equation}\label{7boundNi}
\mathcal{N}^{(a,b)}_{i}\leq\frac{N(N-i)}{N-i-1}\frac{|b-a|}{\mathcal{L}}e^{2\|A\| T}\sum_{m=i+1}^{N-1}\|A^{m+1}b\| + N-i-1.
\end{equation}
We define
\[
I_{i}^{(a,b)}=\lbrace{s\in (a,b)\ |\ |g^{(i)(s)}|\geq c(s)\rbrace}
\]
and observe that $I^{(a,b)}_{i}$ is the union of at most $\mathcal{N}_{i}^{(a,b)}+1$ intervals.\\
We finally set $I_{i}$ to be the union of the $I_{i}^{(a_j,b_j)}$ over all the (at most $\mathcal{N}_{i-1}$)
connected components $(a_j,b_j)$ of $J_{i-1}$.
Therefore,  $I_{i}$ is the union of at most $\mathcal{N}_{i}+\mathcal{N}_{i-1}$ intervals, where
\begin{equation}\label{Ni+1}
\mathcal{N}_{i}=\sum_{j=1}^{\mathcal{N}_{i-1}}\mathcal{N}_{i}^{(a_j,b_j)}\leq
\frac{N(N-i)}{N-i-1}\frac{T}{\mathcal{L}}e^{2\|A\|T}\sum_{m=i+1}^{N-1}\|A^{m+1}b\| +  \mathcal{N}_{i-1}  (N-i-1). 
\end{equation}
Finally we observe that $J_{i}$ is the union of at most $\mathcal{N}_{i}$ open intervals.\\
\indent If $i=N-1$, we observe that for each $s\in J_{N-2}$, recalling (ind2) we have $|g^{(N-1)}(s)|\geq c(s)$.
Therefore we set $J_{N-1}=\emptyset$ and $I_{N-1}=J_{N-2}$. The proof is concluded.
\qed

\bigskip

\noindent We are now going to prove the main result of this subsection.
\begin{Theorem}\label{7LA}
Consider the linear control system
\begin{equation}\label{7linsys}
\dot{x}=Ax+Bu,
\end{equation}
where $1\leq M\leq N$, $A\in\mathbb{M}_{N\times N}$, $B\in\mathbb{M}_{N\times M}$, and $u=(u_1,\ldots ,u_M)\in [-1,1]^M$.\\
Assume that (\ref{7linsys}) is normal, i.e., for every column $b_j$, $j=1,\ldots ,M$, of $B$,
\[
\mathrm{rank}\ [b_j,Ab_j,\ldots ,A^{N-1}b_j]=N.
\]
Let $\mathcal{R}^T$ be defined according to \eqref{reprRT}.
Then for all $T>0$ there exists a constant $\gamma>0$, depending only on $N,M,A,B,T$ such that for all $x,y\in\mathcal{R}^T$,
for all $\zeta\in N_{\mathcal{R}^T}(x)$, the inequality
\begin{equation}\label{7power}
\langle\zeta,y-x\rangle\leq -\gamma\ \|\zeta\|\ \|y-x\|^N
\end{equation}
holds. Moreover, there exists another constant $\gamma'$, depending only on $N,M,A,B,T$, such that
\begin{equation} \label{ballaut}
\text{the ball }\ B(0,\gamma' T^N)\ \text{ is contained in } \mathcal{R}^{T}.
\end{equation}
for all $T>0$. Finally, the constants $\gamma$ and $\gamma'$ are bounded away from zero as $T\to 0^+$. 
\end{Theorem}
\textbf{Proof.}
We consider first the case $M=1$, so (\ref{7linsys}) reads as
\[
\dot{x}=Ax+bu\quad,\ |u|\leq 1,
\]
for a suitable $b\in\R^N$.
Fix $\bar{x}\in\mathrm{bdry}\,\mathcal{R}^T$ together with an optimal control $\bar{u}(\cdot)$ steering $0$ to $\bar{x}$ in time $T$
and $\zeta\in N_{\mathcal{R}^T}(\bar{x})$, $\|\zeta \|=1$. We assume first that $\zeta$ satisfies
Pontryagin's maximum principle, i.e., for a.e. $t\in [0,T]$,
\[
\bar{u}(t)=\mathrm{sign}\langle\zeta,e^{A(T-t)}b\rangle
\]
and prove \eqref{7power} for such $\zeta$. The general case will then follow using Proposition \ref{7convex}.\\
So, we take $\bar{y}\in\mathcal{R}^T$ together with a control $u(\cdot)$ steering $0$ to $\bar{y}$ and compute:
\begin{eqnarray*}
\langle\zeta,\bar{y}-\bar{x}\rangle &=&\int_0^T\langle\zeta, e^{A(T-t)}b\rangle(u(t)-\bar{u}(t))\, dt\\
&&\qquad \text{(recalling \eqref{signlin})}\\
&=&-\int_0^T |\langle \zeta, e^{A(T-t)}b\rangle|\, |u(t)-\bar{u}(t) |\, dt.
\end{eqnarray*}
Set $K(t)=\frac{1}{2}|u(t)-\bar{u}(t)|$ and observe that $0\leq K(t)\leq 1$ for a.e. $t\in [0,T]$, and 
\begin{equation}\label{7ls}
\langle\zeta,\bar{y}-\bar{x}\rangle =-2\int_{0}^{T}|\langle\zeta,e^{A(T-t)}b\rangle|K(t)dt=-2\int_0^T|\langle\zeta,e^{At}b\rangle| K(T-t)dt.
\end{equation}
Moreover, 
\begin{equation}\label{7rs}
\|\bar{y}-\bar{x}\|=\left\|\int_0^Te^{A(T-t)}b(u(t)-\bar{u}(t))dt\right\|\leq 2e^{T\|A\|}\|b\|\int_0^TK(t)dt.
\end{equation}
Set, for $s\in [0,+\infty)$,
\[
g(s)=\langle e^{As}b\ ,\ \zeta\rangle.
\]
By Lemma \ref{7fin int} there exist disjoint sets $I_0,I_1,\ldots ,I_{N-1}$ and numbers $\mathcal{N}_i$ such that $[0,T]=\bigcup_{i=0}^{N-1}I_i$,
each $I_i$ is the disjoint union of at most $\mathcal{N}_i$ intervals and (\ref{7estgi}) holds. Observe that, in particular, it follows
that $g$ may vanish at most at finitely many times, and so, recalling \eqref{signlin}, the control $\bar{u}$ is
piecewise constant and equal to either $1$ or to $-1$.\\
We rewrite 
\begin{equation}\label{7lssum}
\langle\zeta , \bar{y}-\bar{x}\rangle =-2\sum_{i=0}^{N-1}\int_{I_i}|g(s)|K_1(s)ds
\end{equation}
where $K_1(s)=K(T-s)$. We are now going to estimate separately the integrals $\int_{I_i}|g(s)|K_1(s)ds$, for all $i=0,1,\ldots ,N-1.$\\
For $i=0$, we have
\begin{equation}\label{7est i=0}
\int_{I_0}|g(s)|K_1(s)ds\geq\frac{\mathcal{L}}{N}e^{-\|A\|T}\int_{I_0}K_1(s)ds.
\end{equation}
Fix now $i=1,2,\ldots ,N-1$, and write, recalling Lemma \ref{7fin int}, 
\[
\overline{I_i}=\bigcup_{j=1}^{\mathcal{N}_i}[a_{ij},b_{ij}]
\]
where all \textit{open} intervals $(a_{ij},b_{ij})$ are disjoint. Recalling (\ref{7estgi}), we have, for all
$s\in I_i$, $|g^{(i)}(s)|\geq\frac{\mathcal{L}}{N}e^{-\|A\|T}$. Fix $j\in\lbrace{1,2,\ldots ,\mathcal{N}_i\rbrace}$.
We are now going to apply inductively Lemma \ref{7ineqg} on $[a_{ij},b_{ij}]$ with the functions $g^{(i-k-1)}$ in place of $f$,
for $k=0,\ldots ,i-1$. Let $k=0$ and set $f=g^{(i-1)}$. Then the assumption (\ref{7bound a}) is satisfied with $C=\frac{\mathcal{L}}{N}e^{-\|A\|T}$,
thanks to \eqref{7estgi}, and Lemma \ref{7ineqg} yields that for some point $c^0_{ij}\in [a_{ij},b_{ij}]$ we have 
\[
|g^{(i-1)}(s)|\geq C|s-c^0_{ij}|\quad\forall s\in [a_{ij},b_{ij}].
\]
Let $k=1$. By applying Lemma \ref{7ineqg} on each of the two (possibly degenerate) intervals $a_{ij},c^0_{ij}$, $c^0_{ij},b_{ij}$
to the function $f=g^{(i-2)}$ with $C=\frac{\mathcal{L}}{N}e^{-\|A\|T}$, we find suitable points $c^1_{ij}\in [a_{ij},c^0_{ij}]$
and $c^2_{ij}\in [c^0_{ij},b_{ij}]$ such that we have both
\begin{eqnarray*}
|g^{(i-2)}(s)|&\geq &\frac{C}{2}(s-c^1_{ij})^2\quad\forall s\in [a_{ij},c^0_{ij}]
\end{eqnarray*}
and
\begin{eqnarray*}
|g^{(i-2)}(s)|&\geq &\frac{C}{2}(s-c^2_{ij})^2\quad\forall s\in [c^0_{ij},b_{ij}].
\end{eqnarray*}
By continuing the induction process until $k=i-1$, we split the interval $[a_{ij},b_{ij}]$ into at most $2^{i}$ intervals
$[a_{ij}=c^0_{ij},c^1_{ij}],[c^1_{ij},c^2_{ij}],\ldots ,[c^{2^i-1}_{ij},b_{ij}=c^{2^i}_{ij}]$ (some of them being possibly
degenerate) such that for all $l=0,1,\ldots ,2^i-1$ and $s\in [c^l_{ij},c^{l+1}_{ij}]$ one has
\begin{equation}\label{7estij}
\mathrm{either}\ \ |g(s)|\geq\frac{C}{i!}(s-c^l_{ij})^i\quad\mathrm{or}\ \ |g(s)|\geq\frac{C}{i!}(c^{l+1}_{ij}-s)^i.
\end{equation}
Recalling (\ref{7lssum}) and the above discussion, we have 
\begin{eqnarray*}
\langle\zeta,\bar{y}-\bar{x}\rangle &=&-2\sum_{i=0}^{N-1}\int_{I_i}|g(s)|K_1(s)\, ds\\
&=&-2\Bigg[\int_{I_0}|g(s)|K_1(s)\, ds+\sum_{i=1}^{N-1}\sum_{j=0}^{\mathcal{N}_i-1}\sum_{l=0}^{2^i-1}\int_{c^l_{ij}}^{c^{l+1}_{ij}}|g(s)|K_1(s)\, ds\Bigg].
\end{eqnarray*}
Recalling (\ref{7est i=0}) and (\ref{7estij}), and $C=\frac{\mathcal{L}}{N}e^{-\|A\|T}$, we obtain from the above expression that
\begin{equation}\label{7bigsum}
\langle\zeta,\bar{y}-\bar{x}\rangle\leq -\frac{2\mathcal{L}}{N}e^{-\|A\| T}\Bigg[\int_{I_0}K_1(s)\, ds+\sum_{i=1}^{N-1}
\sum_{j=0}^{\mathcal{N}_i-1}\sum_{l=0}^{2^i-1}\int_{c^l_{ij}}^{c^{l+1}_{ij}}\frac{|s-\bar{c}^l_{ij}|^i}{i!}K_1(s)\, ds\Bigg]
\end{equation}
where $\bar{c}^l_{ij}$ is either $c^{l}_{ij}$ or $c^{l+1}_{ij}$, according to the two possibilities appearing in (\ref{7estij}).
Applying Lemma \ref{7app-intineq} to each summand of (\ref{7bigsum}) we therefore obtain
\begin{eqnarray*}
\langle\zeta,\bar{y}-\bar{x}\rangle &\leq & -2\frac{\mathcal{L}}{N}e^{-\|A\| T}\Bigg[\int_{I_0}K_1(s)\, ds+
  \sum_{i=1}^{N-1}\sum_{j=0}^{\mathcal{N}_i-1}\sum_{l=0}^{2^i-1}\frac{\big(\int_{c^l_{ij}}^{c^{l+1}_{ij}}K_1(s)\, ds\big)^{i+1}}{(i+1)!} \Bigg]\\
& &\text{(using the convexity of $x\mapsto x^{i+1}$ on the positive half line)}\\
&\leq &-\frac{2\mathcal{L}}{N}e^{-\|A\| T}\Bigg[\int_{I_0}K_1(s)\, ds+\sum_{i=1}^{N-1}\sum_{j=0}^{\mathcal{N}_i-1}\frac{1}{(i+1)!\ 2^{i^2}}
  \Big(\int_{a_{ij}}^{b_{ij}}K_1(s)\, ds\Big)^{i+1}\Bigg]\\
&\leq & -\frac{2\mathcal{L}}{N}e^{-\|A\| T} \Bigg[\int_{I_0}K_1(s)\, ds+\sum_{i=1}^{N-1}\frac{1}{(i+1)!\ 2^{i^2}\mathcal{N}_i^i}
\Big(\int_{I_i}K_1(s)\, ds\Big)^{i+1}\Bigg].
\end{eqnarray*}
Thus, recalling that $0\le K_1(s) \le 1$ a.e.,
\begin{eqnarray}\label{7last ineq}
\begin{aligned}
\langle\zeta,\bar{y}-\bar{x}\rangle\; & \leq \;  -\frac{2\mathcal{L}}{N}e^{-\|A\| T} \Bigg[\frac{1}{|I_0|^{N-1}}\Big(\int_{I_0}K_1(s)\, ds\Big)^N \\
&\qquad  +\sum_{i=1}^{N-1}\frac{1}{(i+1)!\ 2^{i^2}\mathcal{N}_i^i|I_i|^{N-i-1}}\Big(\int_{I_i}K_1(s)\, ds\Big)^N\Bigg].
\end{aligned}
\end{eqnarray}
Owing to (\ref{7boundN0}) and (\ref{Ni+1}), we see that $\mathcal{N}_0 < \mathcal{N}_1 < \ldots < \mathcal{N}_N\leq C(A,b,N)e^{2\|A\|T}$
where $C(A,b,N)$ depends only on $A,b,N$. Therefore, we obtain finally from (\ref{7last ineq}) and the definition of $K_1$ that 
\begin{equation}\label{7lsrs}
\langle\zeta,\bar{y}-\bar{x}\rangle\leq -C(A,b,N,T)e^{-\|A\|T}\Big(\int_{0}^T|u(t)-\bar{u}(t)|\Big)^N,
\end{equation}
where $C(A,b,N,T)$ is a positive constant, depending only on $A,b,N,T$ such that 
\begin{equation}\label{7lboundC}
\liminf_{T\rightarrow 0}C(A,b,N,T)>0.
\end{equation}
Recalling (\ref{7rs}), we complete the proof for the case $M=1$ (i.e., a scalar control) by setting 
\[
\gamma=2^{-N}e^{-\| A\|(N+1)T}\|b\|^{-N}C(A,b,N,T).
\]
\indent Let now $M>1$. Take $\bar{x}\in\mathrm{bdry}\,\mathcal{R}^T$ together with an optimal control
$\bar{u}(\cdot)=(\bar{u}_1(\cdot),\ldots ,\bar{u}_{M}(\cdot))$ steering the origin to $\bar{x}$ in the optimal time $T$,
and $\bar{y}\in\mathcal{R}^T$ with a control $u(\cdot)=(u_1(\cdot),\ldots ,u_M(\cdot))$ steering the origin to $\bar{y}$ in time $T$.
Then, for each $\zeta\in N_{\mathcal{R}^T}(\bar{x})$, $\|\zeta\|=1$, we can write
\begin{equation}\label{7lsm1}
\langle\zeta,\bar{y}-\bar{x}\rangle \leq \int_0^T\big\langle \zeta, e^{A(T-s)}Bw(s)\big\rangle ds =\sum_{i=1}^{M}\int_0^T
   \big\langle\zeta , e^{A(T-s)}b_i\big\rangle (u_i(s)-\bar{u}_i(s))\, ds,
\end{equation}
where $B=\big(b_i\big)_{i=1,2,\ldots ,M}$. Recalling \eqref{signlin} we have also 
\begin{equation}\label{7lsm2}
\langle\zeta,\bar{y}-\bar{x}\rangle =-\sum_{i=1}^M\int_0^T\big|\big\langle\zeta,e^{A(T-s)b_i}\big\rangle\big|\, \big|u_i(s)-\bar{u}_i(s)\big|\, ds.
\end{equation}
Moreover,
\begin{equation}\label{7rsm}
\|\bar{y}-\bar{x}\|\leq e^{\|A\|T}\sum_{i=1}^M\|b_i\|\int_0^T|u_i(s)-\bar{u}_i(s)|ds.
\end{equation}
We now apply the same argument leading to (\ref{7lsrs}) to each summand of the right hand side of (\ref{7lsm1}).
Therefore we obtain, using (\ref{7lsm2}), that 
\begin{equation}  \label{eq:forball}
\langle\zeta,\bar{y}-\bar{x}\rangle\leq -C'(A,B,T,N,M)e^{-||A\| T}\sum_{i=1}^M\Big(\int_0^T|u_i(s)-\bar{u}_i(s)|\, ds\Big)^N,
\end{equation}
where the positive constant $C'$ depends only on $A,B,T,N,M$ and 
\[
\liminf_{T\rightarrow 0}C'(A,B,T,N,M)>0.
\]
We conclude the proof of \eqref{7power} by applying (\ref{7rsm}) and setting 
\[
\gamma=2^Ne^{-(N+1)\|A\|T} C''(A,B,T,N,M),
\]
where $C''$ is a constant enjoying the same properties as $C'$.\\
\indent In order to prove the statement concerning the ball contained in $\mathcal{R}^T$, observe that the inequality \eqref{eq:forball}
with $u\equiv 0$ and $\bar{y}=0$, taking into account that the control $\bar{u}$ is necessarily bang-bang, becomes
\[
\langle \zeta ,\bar{x}\rangle \ge C'(A,B,T,N,M) e^{-\|A\|T}MT^N, 
\]
from which (taking $\| \zeta\|=1$) we obtain
\begin{equation} \label{0notin}
 \| \bar{x}\| \ge C'(A,B,T,N,M) e^{-\|A\|T}MT^N:=\gamma' T^N.
\end{equation}
The above inequality yields in particular that $0$ belongs to the interior of $\mathcal{R}^T$.
Since $\mathcal{R}^T$ is convex and \eqref{0notin} holds for all $\bar{x}\in\mathrm{bdry}\,\mathcal{R}^T$, \eqref{ballaut} follows.\\
The last statement follows from \eqref{7lboundC} and the explicit expressions for $\gamma$ and $\gamma'$. The proof is concluded.
\qed
\begin{Rem}
The exponent $N$ in (\ref{7power}) is optimal.
\end{Rem}
In fact, consider the dynamics
\[
x^{(N)}=u,\qquad u\in [-1,1]. 
\]
Let $x_1(\cdot)$ be the solution corresponding to the control $u\equiv 1$. Fix $s>0$ and let $x_s(\cdot)$ be the solution corresponding to the control
\[
u_s(t)=
\begin{cases}
 1 & 0<t<s\\
-1& s<t.
\end{cases}
\]
Fix any $T>0$ and observe that $x_1(T) = T^N/N!$, while, for $0<s<T$,
\[
x_s(T)-x_1(T) = \frac{-2}{N!}(T-s)^N. 
\]
Observe also that $|x_s^{(N-1)}(T)-x_1^{(N-1)}(T)|=2|T-s|$, so that setting $X_i=(x_i,x_i',\ldots ,x_i^{(N-1)})$, $i=1,s$, one obtains
\[
\langle e_1, X_s(T)-X_1(T)\rangle = \frac{-2}{N!} (T-s)^N \ge - \gamma \big\| X_s(T)-X_1(T)  \big\|^N, 
\]
for a suitable positive constant $\gamma$. If $s\to T$, then $X_s(T)\rightarrow X_1(T)$, and this shows that $N$ is the smallest exponent allowed in
\eqref{7power}.\qed
\subsection{Nonautonomous systems}\quad\\
\indent The following Lemma is a first step for studying the reachable sets in the case of nonlinear control systems by using
the linearization approach which we design in this paper. We will prove that under the rank condition (normality) at 0 of
the linear nonautonomous control system (\ref{7linmo}),
the strict convexity of the reachable sets is preserved up to a sufficiently small time, provided $A(t)$ and $B(t)$ are not too
far from $A(0)$, $B(0)$. 
Sufficient conditions for the validity of the assumptions if $N=2$ will be given below. In Section \ref{sec:bangbang} sufficient
conditions in the case where $A$ and $B$ come from a linearization around a trajectory will also be given.

Let $A:\R^+\rightarrow \mathbb{M}_{N\times N}$ and $B:\R^+\rightarrow \mathbb{M}_{N\times M}$, $1\le M\le N$, 
be measurable and consider the linear nonautonomous control system
\begin{equation}\label{7linmo}
\dot{x}(t)=A(t)x(t)+B(t)u(t),
\end{equation}
together with
\begin{equation}\label{7linm0}
\dot{x}(t)=A(0)x(t)+B(0)u(t),
\end{equation}
where $u=(u_1,\ldots ,u_M)\in [-1,1]^M$. We denote by $M^\top(\cdot,\cdot)$ the matrix solution of
\begin{equation}\label{7defM}
\left\{\begin{array}{ll}
\frac{\partial}{\partial t}X(t,s)=X(t,s)A(t)\ \mathrm{for}\ t,s\geq0\\
M(s,s)=\mathbb{I}
\end{array}\right.
\end{equation}
and by $M_0^\top(\cdot,\cdot)$ the matrix solution of
\begin{equation}\label{7defM0}
\left\{\begin{array}{ll}
\frac{\partial}{\partial t}X(t,s)=X(t,s)A(0)\ \mathrm{for}\ t,s\geq0\\
M_0(s,s)=\mathbb{I}.
\end{array}\right.
\end{equation}
Let $b(\cdot)$ be a column of $B(\cdot)$, let $\mathcal{T}>0$ and define, for any $\zeta\in\mathbb{R}^N$,
\begin{equation}\label{7defgg}
g(t)=\langle \zeta M^\top(\mathcal{T},t),b(t)\rangle,
\end{equation}
and 
\begin{equation}\label{7defg0}
g_0(t)=\langle \zeta M_0^\top(\mathcal{T},t),b(0)\rangle.
\end{equation}
Observe that $g$ and $g_0$ can be seen as the switching functions related to \eqref{7linmo} and \eqref{7linm0}, respectively.\\
We state now an abstract result which permits to transfer to $g$ some properties of $g_0$ and to establish the quantitative strict
convexity estimate for the reachable set from the origin of \eqref{7linmo}. Sufficient conditions in order to apply the following lemma
to suitable linearizations in dimension $2$ and $3$ will be given in Section \ref{sec:bangbang}.
\begin{Lemma}\label{7gg0}
Fix $\zeta\in\mathbb{R}^N$, $\|\zeta\|=1$, and let $g$ and $g_0$ be defined according to (\ref{7defgg}) and \eqref{7defg0}, respectively.
Assume that
\begin{itemize}
\item[(i)] $\mathrm{rank}[b(0),A(0)b(0),\ldots,A^{N-1}(0)b(0)]=N $;
\item[(ii)] $g$ is $N-2$ times differentiable and $g^{(N-2)}$ is absolutely continuous;
\item[(iii)] there exists a constant $K=K(A,B)$ and a time $\mathcal{T}'$ such that for all $i=0,\ldots ,N-1$ one has
\[
 |g^{(i)} (t)-g_0^{(i)}(t)|\le K t\qquad\text{ for all  }\; 0\le t\le\mathcal{T}'.
\]
\end{itemize}
Let $\mathcal{L}$ be defined by \eqref{7lbound}, with $A(0)$, $b(0)$ in place of $A$, $b$, respectively.\\
Then there exists $\mathcal{T}=\mathcal{T}(A,B,N,\mathcal{T}')$ with the following property:\\
for every $0 < \tau<\mathcal{T}$ there exist disjoint sets $I_i$ and numbers $\mathcal{N}_i$, $i=0,\ldots ,N-1$, depending only on $A(0),B(0),\tau$, such that 
\begin{enumerate}
\item [(a)] $[0,\mathcal{T}]=\bigcup_{i=0}^{N-1}I_i$;
\item [(b)] each $I_i$, $i=0,\ldots,N-1$ is the disjoint union of at most $\mathcal{N}_i$ intervals;
\item [(c)] for each $i=0,\ldots,N-2$ and each $s\in I_i$
\[
|g^{(i)}(s)|\geq\frac{\mathcal{L}}{2N}e^{-\|A(0)\|s};
\]
\item [(d)] $g^{(N-1)}$ has constant sign in every connected component of $I_{N-1}$ and, for each $s\in I_{N-1}$
\[
|g^{(N-1)}(s)|\geq \frac{\mathcal{L}}{2N}e^{-\|A(0)\|s}.
\]
\end{enumerate}
\end{Lemma}
\begin{Rem}
$I_0,\ldots ,I_{N-1}$ are exactly the intervals constructed in Lemma \ref{7fin int} with $g_0$ in place of $g$.
\end{Rem}
\noindent \textbf{Proof of Lemma \ref{7gg0}.} 
Let $I_0,\ldots,I_{N-1}$ be the sets appearing in the statement of Lemma \ref{7fin int}, with $g_0$ in place of $g$, so that, in particular, (a) and (b) hold.
Moreover, for each $s\in I_i$ we have
\[
 |g_0^{(i)} (s)| \ge \frac{\mathcal{L}}{N}e^{-\| A(0)\| \tau}.
\]
Therefore, if $K\mathcal{T} < \frac{\mathcal{L}}{2N}e^{-\| A(0)\| s}$, also (c) and (d) hold, owing to (iii).
\qed
\begin{Rem}\label{remsuffcond}
Sufficient conditions for the validity of the assumptions of Lemma \ref{7gg0} in the case $N=2$.
\end{Rem}
Let $N=2$. Sufficient conditions for the validity of assumptions (ii) and (iii) in the above Lemma are the following:
\begin{itemize}
 \item[($C_0$)] $A(\cdot)$ is measurable and 
\[
\|A(t)-A(0)\|\leq Lt\ \mathrm{for\ all}\ t\geq 0;
\]
\item[($C_1$)] $B(\cdot)$ is absolutely continuous and 
\[
\Big\|\frac{d}{dt}B(t)\Big\|\leq 2Lt\ \mathrm{for\ all}\ t\geq 0,
\]
where $L$ is a positive constant.
\end{itemize}
Observe that condition $(C_0)$ implies that $t=0$ is a continuity point for $A(\cdot)$, so that $A(0)$ in $(C_0)$ is meaningful.\qed

\medskip

\noindent As an immediate corollary of Lemma \ref{7gg0} we obtain the following
\begin{Theorem}\label{7key-nonauton} Consider the linear nonautonomous control system \eqref{7linmo} under the assumptions of Lemma
\eqref{7gg0}.
Let $\mathcal{R}^{\tau}$ denote the reachable set at time $\tau >0$ from the origin for \eqref{7linmo}.
Then there exist a time $\mathcal{T}=\mathcal{T}(A,B,N)>0$ and a constant $\gamma=\gamma(A,B,N)>0$ such
that for every $0\leq\tau\leq\mathcal{T}$, for every $x,y\in\mathcal{R}^{\tau}$, for every $\zeta\in N_{\mathcal{R}^{\tau}}(x)$, we have
\[
\langle\zeta,y-x\rangle\leq -\gamma \|\zeta\|\|y-x\|^N.
\]
Moreover, there exists another constant $\gamma'=\gamma'(A,B,N)>0 $ such that for every $0 < \tau \le \mathcal{T}$ 
\begin{equation} \label{ballnonaut}
\text{the ball }\ B(0,\gamma' \tau^N)\ \text{ is contained in } \mathcal{R}^{\tau}.
\end{equation}
\end{Theorem}

\medskip

\noindent \textbf{Proof.} The argument developed in the proof of Theorem \ref{7LA} can be used also in this case.
Indeed, fix $\bar{x}\in\mathrm{bdry}\,\mathcal{R}^{\tau}$ together with a control $\bar{u}(\cdot)$ steering $0$ to
$\bar{x}$ in time $\tau$ and let $\zeta\in N_{\mathcal{R}^{\tau}}(\bar{x})$, $\|\zeta\|=1$, be such that
\begin{equation}\label{usign0}
 \bar{u}(t) = \mathrm{sign}\langle \zeta , b(t)M(\tau,t)\rangle\qquad \text{for a.e. }\; t\in (0,\tau)
\end{equation}
(here, as at the beginning of the proof of Theorem \ref{7LA}, we assume that $B=b$ is a vector, i.e., the control is scalar).\\
Then the proof proceeds exactly as for Theorem \ref{7LA}, provided that $g$ is given by
\[
 g(t)=\langle \zeta , b(t)M(\tau,t)\rangle,
\]
so that for all $\bar{y}\in\mathcal{R}^{\tau}$ one has
\[
 \langle \zeta , \bar{y}-\bar{x}\rangle = - 2 \sum_{i=0}^{N-1} \int_{I_i} |g(s)|K_1(s)\, ds,
\]
where $K_1(s)=\frac12 |u(\tau -s)-\bar{u}(\tau-s)|$, $u(\cdot)$ is the control which steers the origin to $\bar{y}$, and the
sets $I_i$ are those appearing in the statement of Lemma \ref{7gg0}.\qed


%
\section{A nonlinear bang bang principle in dimensions $2$ and $3$}\label{sec:bangbang}
Starting from the present section we will deal with nonlinear control systems, which are affine and symmetric with respect to the control.
This section is devoted to giving sufficient conditions so that controls steering the origin to the boundary of the reachable
set are always bang-bang, provided the final time is sufficiently small. More precisely, the following result holds.
\begin{Theorem}\label{bangbang}
Let $N=2$ or $N=3$. Consider the control system
\begin{equation} \label{NLsys}
\begin{cases}
 \dot{x}(t)=&F(x(t)) + G(x(t)) u(t),\\
x(0)=&0,
\end{cases}
\end{equation}
where $1\le M\le N$, $u(\cdot)=(u_1(\cdot),\ldots,u_M(\cdot))\in[-1,1]^M$ a.e., and $F:\mathbb{R}^N\rightarrow\mathbb{R}^N$ and
$G=(G_1,\ldots,G_M):\mathbb{R}^N\rightarrow\mathbb{M}_{N\times M}$ satisfy the following assumptions:
\begin{itemize}
\item[(i)] $F$ and $G$ are of class $\mathcal{C}^{N-1}$ and all partial derivatives are Lipschitz with constant $L$;
\item[(ii)] $F(0)=0$;
\item[(iii)] $\mathrm{rank}[G_i(0),DF(0)G_i(0),\ldots,(DF(0))^{N-1}G_i(0)]=N$\, for\, $i=1,\ldots,M$;
\item[(iv)] if $N=3$, then $DG(0)=0$ and $D^2F(0)=0$.
\end{itemize}
Let $\mathcal{R}^\tau$ denote the reachable set of \eqref{NLsys} at time $\tau >0$. Then there exists $\mathcal{T}>0$, depending only on $DF(0),G(0),L,N$,
such that for every $0<\tau<\mathcal{T}$ the following properties hold:
\begin{itemize}
\item[(a)] every admissible control $u(\cdot)$ such that the corresponding trajectory $y^u(\cdot)$ of \eqref{NLsys} at time $\tau$ belongs to the boundary of
$\mathcal{R}^\tau$ is essentially determined by the curve $\lambda(\cdot)$, the solution of the adjoint equation
\begin{equation}\label{adjeq}
 \begin{cases}
  \dot{\lambda}(t)&=-\lambda(t) \big(DF(y^u(t))+DG(y^u(t))u(t)\big)\\
  \lambda(\tau)&=\zeta\in N_{\mathcal{R}^\tau}(y^u(\tau)),\quad\|\zeta\|=1,
 \end{cases}
\end{equation}
through the identity
\begin{equation}\label{maxu}
 u(t)=\mathrm{sign}\,\langle \lambda (t), G(y^u(t))\rangle\quad \text{ a.e.;}
\end{equation}
\item[(b)] $u$ is bang-bang, i.e., $u(t)\in\{ -1,1\}^M$ a.e., and is piecewise constant;
\item[(c)] the maximum number of discontinuities of $u$ depends only on $DF(0)$, $G(0)$, $L$, and $N$.
\end{itemize}
\end{Theorem}
\textbf{Proof.} We consider first the case where $M=1$, i.e., $G(x)$ is a vector and the control $u$ is scalar.\\
Fix $\tau >0$ and an admissible control $u$ such that $\bar{x}:=y^u(\tau)\in\mathrm{bdry}\,\mathcal{R}^{\tau}$. By the Maximum Principle (see Theorem
\ref{th:PMPnonlin}) there exist $\zeta\in N_{\mathcal{R}^\tau}(\bar{x})$, $\|\zeta\|=1$, and an adjoint curve $\lambda(\cdot)$, a solution of
\eqref{adjeq}, such that $u$ satisfies \eqref{maxu}. Proving (a), (b), and (c) amounts to showing that the switching function $\langle\lambda(t),G(y^u(t))\rangle$
vanishes at most finitely many times in $[0,\tau]$ and the number of its zeros depends only on $DF(0)$, $G(0)$, $L$, and $N$.\\
For convenience we rewrite the switching function in the following way. Let $M^\top(\cdot,\cdot)$ denote the matrix solution of \eqref{7defM},
with $A(t)=DF(y^u(t))+DG(y^u(t))u(t)$, and set $b(t)=G^\top(y^u(t))$, $t\in [0,\tau]$. Then
\[
g(t):= \langle\lambda (t),b(t)\rangle = \langle\zeta, b(t)M(\tau,t)\rangle. 
\]
Let also $M_0^\top(\cdot,\cdot)$ be defined by \eqref{7defM0}, where $A(0)=DF(0)$, and let $g_0$ be defined according to \eqref{7defg0}.
We wish to apply Lemma \ref{7gg0} to the above introduced mappings $g$ and $g_0$.\\
First of all, we compute $g'$ and observe that it is continuous. Indeed,
\begin{eqnarray*}
 g'(t)&=&\big\langle\zeta,b'(t)M(\tau,t)\rangle + \langle \zeta,b(t)\frac{\partial}{\partial t}M(\tau,t)\big\rangle\\
      &=&\big\langle\zeta, DG^\top(y^u(t))\dot{y}^u(t)M(\tau,t)\big\rangle + \big\langle\zeta,-b(t)A^\top(t)M(\tau,t)\rangle\\
      &=& \big\langle\zeta,DG^\top(y^u(t))\big( F(y^u(t))+G(y^u(t))u(t)\big)^\top M(\tau,t)\big\rangle\\ 
      && \qquad + \big\langle\zeta,-G^\top (y^u(t))\big( DF(y^u(t))+DG(y^u(t))u(t)\big)^\top M(\tau,t)\big\rangle\\
      &=&\big\langle\zeta ,\big(DG^\top(y^u(t))F^\top(y^u(t))-DF^\top(y^u(t))G^\top(y^u(t))\big)M(\tau,t)\big\rangle\\
      &=&\big\langle \lambda(t),[F^\top,G^\top](y^u(t))\big\rangle,
\end{eqnarray*}
where $[F^\top,G^\top](x):=DG^\top(x)F^\top(x)-DF^\top(x)G^\top(x)$ denotes the Lie bracket.\\
Moreover, if $N=3$ $g'$ is a.e. differentiable and we have
\[
 g''(t) = \big\langle \dot{\lambda}(t),[F^\top,G^\top](y^u(t))\big\rangle + \big\langle \lambda (t),\frac{d}{dt} [F^\top,G^\top](y^u(t))\big\rangle.
\]
Finally, observing that $g_0'(t)=\big\langle\zeta , -DF^\top (0)G^\top (0)M_0(\tau,t)\big\rangle=
\big\langle\zeta M_0^\top(\tau,t),[F^\top,G^\top](0)\big\rangle$, we have
\begin{eqnarray*}
g'(t)-g_0'(t)&=&\big\langle\zeta \big( M^\top (\tau,t)-M_0^\top(\tau,t)\big),[F^\top,G^\top](y^u(t))\big\rangle\\
&&\qquad + \big\langle\zeta M_0^\top(\tau,t),DG^\top(y^u(t)) F^\top(y^u(t))\big\rangle\\
&&\qquad + \big\langle\zeta, \big( DF^\top(0)-DF^\top(y^u(t))\big)G^\top(y^u(t))\big\rangle\\
&&\qquad + \big\langle\zeta ,DF^\top(0)\big( G^\top(0)-G^\top(y^u(t))\big)\big\rangle\\
&=& I + I\!I + I\!I\!I + I\!V.
\end{eqnarray*}
We are now going to estimate separately each summand of the above expression.\\
First, we observe that, thanks to the assumptions (i) and (ii), we have
\begin{eqnarray*}
 \| M(\tau,t)-M_0(\tau,t)\| &\le & \int_0^t \Big( \| A(s)\| \, \| M_0(\tau,s)-M(\tau,s)\| + \| M_0(\tau,s)\| \, \| A(0)-A(s)\|\Big)\, ds\\
&\le& K_1 \int_0^t \| M(\tau,s)-M_0(\tau,s)\|\, ds + K_2 t^2,
\end{eqnarray*}
where $K_1$ and $K_2$ are suitable constants depending only on $DF(0)$, $G(0)$, $L$, and $\tau$. Gronwall's lemma therefore yields
\[
\| M(\tau,t)-M_0(\tau,t)\| \le K_3t^2\qquad\text{ for all }\; t\in [0,\tau],
\]
where the constant $K_3$ depends only on $K_1$, $K_2$. Therefore, there exists a constant $K_I$, depending only on $DF(0)$, $G(0)$, $L$, and $\tau$ such that
\begin{equation}\label{estI}
|I|\le K_{I} t^2 \qquad\text{ for all }\; t\in [0,\tau].
\end{equation}
Assumptions (i) and (ii) yield in turn
\begin{equation}\label{estII}
| I\!I|\le K_{I\!I} t\qquad\text{ for all }\; t\in [0,\tau], 
\end{equation}
for a suitable constant $K_{II}$ depending only on $DF(0)$, $L$, and $\tau$, and
\begin{equation}\label{estIII}
| I\!I\!I | \le K_{I\!I\!I} t  \qquad\text{ for all }\; t\in [0,\tau],
\end{equation}
\begin{equation}\label{estIV}
|I\!V|\le K_{I\!V} t \qquad\text{ for all }\; t\in [0,\tau], 
\end{equation}
where again $K_{I\!I\!I}$ and $K_{I\!V}$ depend only on $DF(0)$, $G(0)$, $L$, and $\tau$.\\
Therefore, summing \eqref{estI}, \eqref{estII}, \eqref{estIII}, and \eqref{estIV},
we obtain that there exist $K$ and $\mathcal{T}'>0$, depending only on $DF(0)$, $G(0)$, and $L$, such that
\begin{equation}\label{estg'}
| g'(t)-g_0'(t)|\le K t\qquad\text{ for all }\; 0\le t \le \mathcal{T}'. 
\end{equation}
Let now $N=3$, and observe that owing to assumption (iv) each summand $I$, $I\!I$, $I\!I\!I$, $I\!V$, divided by $t^2$, is bounded and a.e.
differentiable, so that
\[
|g''(t)-g_0''(t)|\le K' t  \qquad\text{ for a.e. }\; t\in [0,\mathcal{T}'],
\]
where the constant $K'$ depends only on $DF(0)$, $G(0)$, and $L$ and so does $\mathcal{T}'$.\\
Observe that all the above constants do not depend on $\zeta$.
Therefore, invoking assumption (iii), we can apply Lemma \ref{7gg0} (for $N=2,3$), thus obtaining that there exists
$\mathcal{T} >0$, depending
only on $G(0)$, $DF(0)$, $L$, and $N$ with all the properties (a), (b), (c), and (d), which are exactly those required to
complete the proof.\\
In the general case (i.e., $1<M\le N$), it suffices to apply the above argument to each column of $G$. The proof is concluded.\qed
\begin{Rem}\label{dimrestr1}
On the assumption $N=2$ or $N=3$. 
\end{Rem}
The restriction $N\le 3$ depends on our method for comparing the switching function $g$ for the nonautonomous system coming from the linearization
along an optimal trajectory, with the switching function for autonomous system obtained by linearizing at the origin. This comparison
requires higher order derivatives of $g$, whose existence we are not able to insure if $N>3$.\qed
\section{Quantitative strict convexity of reachable sets and uniqueness of optimal controls: the nonlinear two dimensional case}\label{sec:nonlinstrict}
\indent In this subsection, we will show that, provided the linearization at $0$ satisfies the normality condition and the nonlinear part
is smooth and small enough, the reachable set is strictly convex up to a sufficiently small time. 
\begin{Theorem}\label{7NL}
Consider the control system
\begin{equation}\label{7system4}
\left\{\begin{array}{ll}
\dot{x}(t)=F(x(t))+G(x(t))u(t), \\
x(0)=0,
\end{array}\right.
\end{equation}
under the following assumptions (in the following $M$ is either $1$ or $2$):\\
$u(\cdot)=(u_1(\cdot),u_M(\cdot))\in [-1,1]^M$ a.e., $F:\mathbb{R}^2\rightarrow\mathbb{R}^2$, $G:\mathbb{R}^2\rightarrow\mathbb{M}_{2\times M}$
are of class $\mathcal{C}^{1,1}$ (with Lipschitz constant L) and 
\begin{enumerate}
\item [(i)] $F(0)=0$,
\item [(ii)] $\mathrm{Rank}\ [G_i(0),DF(0)G_i(0)]=2$ for $i=1,M$ where $G=(G_1,G_M)$,
\item [(iii)]$DG(0)=0$.
\end{enumerate}
Let $\mathcal{R}^{\tau}$ denote the reachable set at time $\tau >0$ for \eqref{7system4}.
Then there exists $\mathcal{T}>0$, depending only on $L,DF(0),G(0)$, with the following properties:
\begin{itemize}
\item[(a)] for every $\tau \le \mathcal{T}$ and every $x\in\mathrm{bdry}\,\mathcal{R}^\tau$ there exists one and only one admissible control $u$
steering the origin to $x$ in time $\tau$ (and $u$ is bang-bang with finitely many switchings).
\item[(b)] For every $0<\tau<\mathcal{T}$ the reachable
set $\mathcal{R}^{\tau}$ is strictly convex. More precisely, for every $x_1\in\mathrm{bdry}\,\mathcal{R}^{\tau}$ and $x_2\in\mathcal{R}^{\tau}$,
for every $\zeta\in N^{P}_{\mathcal{R}^{\tau}}(x_1)$, one has
\begin{equation}\label{7final-est}
\langle\zeta,x_2-x_1 \rangle\leq -\gamma \|\zeta\|\|x_2-x_1\|^2.
\end{equation}
where $\gamma$ is a positive constant depending only on $L$, $DF(0)$, $G(0)$.
\item[(c)] There exist another time $\mathcal{T}'>\mathcal{T}$, depending only on $L$, $DF(0)$, $G(0)$, such that for every $0<\tau<\mathcal{T}'$
the reachable set $\mathcal{R}^{\tau}$ has positive reach. More precisely, for every $x_1\in\mathrm{bdry}\,\mathcal{R}^{\tau}$ and $x_2\in\mathcal{R}^{\tau}$,
for every $\zeta\in N^{P}_{\mathcal{R}^{\tau}}(x_1)$, one has
\begin{equation}\label{eq:posreach}
\langle\zeta,x_2-x_1 \rangle\leq \gamma' \|\zeta\|\|x_2-x_1\|^2,
\end{equation}
where $\gamma'$ is a nonnegative constant depending only on $L$, $DF(0)$, $G(0)$.
\item[(d)] There exist another positive constant $\gamma''$
and a positive time $\mathcal{T}''\le\mathcal{T}$, depending only on $L$, $DF(0)$, $G(0)$, such that the ball $B(0,\gamma''\tau^2)$
is contained in $\mathcal{R}^{\tau}$ for all $0<\tau<\mathcal{T}''$.
\item[(e)] The minimum time function is continuous
in $\mathcal{R}^{\tau}$, for all $0<\tau<\mathcal{T}''$.
\end{itemize}
\end{Theorem}
\textbf{Proof.} We begin proving the result for $M=1$, i.e., for a scalar control.\\
\indent Fix $\tau>0$ and $x_1\in\mathrm{bdry}\,\mathcal{R}^{\tau}$, together with an optimal control $u_1(\cdot)$ steering $0$ to $x_1$ and
the associate trajectory $x_1(\cdot)$. Take any $x_2\in\mathcal{R}^{\tau}$ together with $u_2(\cdot)$ steering $0$ to $x_2$ and the associate
trajectory $x_2(\cdot)$, and set $x(t)=x_2(t)-x_1(t)$. Then, for a.e. $t\in [0,\tau]$,
\begin{equation}\label{7diffx}
\dot{x}(t)=A_1(t)x(t)+G(x_1(t))w(t),
\end{equation}
where $w(t)=u_2(t)-u_1(t)$ and 
\[
A_1(t)=\int_0^1DF(x_1(t)+\tau x(t))d\tau+ u_2(t)\int_0^1DG(x_1(t)+\tau x(t))d\tau.
\] 
Let $z(\cdot)$ be the solution of the linear system which is defined by linearizing along the optimal trajectory $x_1(\cdot)$:
\begin{equation}\label{7system5}
\left\{\begin{array}{ll}
\dot{z}(t)=A(t)z(t)+G(x_1(t))w(t),\\
z(0)=0,
\end{array}\right.
\end{equation}
where $A(t)=DF(x_1(t))+DG(x_1(t))u_1(t)$.\\
We have 
\begin{eqnarray*}
\frac{d}{dt}\|x(t)-z(t)\|& \leq & \|A_1(t)x(t)-A(t)z(t)\|\\
&\leq & \|A(t)\|\|x(t)-z(t)\|+\|A_1(t)-A(t)\|\|x(t)\|\\
&\leq & L_1 \|x(t)-z(t)\|+L\|x(t)\|^2,
\end{eqnarray*}
where $L_1=\|DF(0)\|+2Le^{2L\tau}$.\\
Thus, by Gronwall's inequality we get
\begin{equation}\label{75dist1}
\|x(t)-z(t)\|\leq e^{L_1t}L\int_{0}^t
\|x(s)\|^2\, ds.
\end{equation}
On the other hand, observing that 
\[
\frac{d}{dt}\|x(t)\|\leq L_2\|x(t)\|+L_3|w(t)|,
\]
(where $L_2=\|DF(0)\|+4Le^{2L\tau}$ and $L_3=|G(0)|+e^{2L\tau}$) we also have
\begin{equation}\label{75dist2}
\|x(t)\|\leq L_3e^{L_2t}\int_0^t|w(s)|ds.
\end{equation}
From (\ref{75dist1}) and (\ref{75dist2}), one obtains
\begin{equation}\label{75dist3}
\|x(t)-z(t)\|\leq L_4 t\Big(\int_{0}^t|w(s)|ds\Big)^2,
\end{equation}
where $L_4=LL_3^2e^{(L_1+2L_2)\tau}$.\\
\indent Since $\bar{x}_1\in\mathrm{bdry}\,\mathcal{R}^\tau$, by Pontryagin's maximum principle there exists an absolutely continuous function
$\lambda:[0,\tau]\rightarrow\mathbb{R}^2$ with the following properties
\[
\dot{\lambda}(t)=-\lambda(t)A(t),\quad \lambda(\tau)=\zeta\in N_{\mathcal{R}^{\tau}}(\bar{x}_1),\; \zeta\neq 0,
\]
\begin{equation}\label{7sign}
u_1(t)=\text{sign}\langle \lambda(t),G(x_1(t))\rangle.
\end{equation}
We set now $b(t)=G(x_1(t))$ and consider the linear nonautonomous control system 
\begin{equation}\label{7system6}
\left\{\begin{array}{ll}
\dot{y}(t)=A(t)y(t)+b(t)u(t),\\
y(0)=0,
\end{array}\right.
\end{equation}
together with the trajectory $y_1(\cdot)$, corresponding to the control $u_1(\cdot)$. 
Observe first that, thanks to \eqref{7sign}, $x_1$ belongs to the boundary of the reachable set at time $\tau$ for
\eqref{7system6}. Observe moreover that $A(\cdot)$ is measurable, and, since both $F$ and $G$ are Lipschitz with constant $L$ and $DG(0)=0$, we have 
\begin{eqnarray*}
\|A(t)-A(0)\|&=& \|DF(x_1(t))+DG(x_1(t))u_1(t)-DF(0)\|\\
&\leq &2L\|x_1(t)\|\leq 2L\| G(0)\| e^{2Lt}t.
\end{eqnarray*}
Finally, $b'(t)=DG(x_1(t))\dot{x}_1(t)$ so that 
\[
\| b'(t)-b(0)\|=\|b'(t)\|\, \leq\, L\|x_1(t)\|\big(2L\|x_1(t)\|+\| G(0)\|\big)\, \leq\, Kt,
\]
where $K=Le^{2L\tau}\big(2Le^{2L\tau}+1\big)\| G(0)\|^2$.\\
Therefore, all assumptions of Theorem \ref{7key-nonauton} are satisfied by \eqref{7system6}, so that there exists a time
$\mathcal{T}_0>0$ depending only on $L,DF(0),G(0)$ such that if $0\leq\tau<\mathcal{T}_0$ the following properties hold:
\begin{itemize}
 \item[i)] $u_1(t)$ is uniquely determined by \eqref{7sign} on $(0,\tau)$,
\item[ii)] there exists a constant $\gamma(\tau)>0$, depending only on $L,DF(0),G(0),\tau$ such that, for all
trajectories $y(\cdot)$ of \eqref{7system6},
$\langle\zeta,y(\tau)-y_1(\tau)\rangle \leq -\gamma(\tau) \|y(\tau)-y_1(\tau)\|^2$. More precisely, recalling (\ref{7lsrs}),
\begin{equation}\label{7linaux}
\langle\zeta,y_2(\tau)-y_1(\tau)\rangle\leq -\gamma_1(\tau) \Big(\int_{0}^{\tau}|w(s)|ds\Big)^2,
\end{equation}
where $y_2(\cdot)$ is the trajectory of (\ref{7system6}) associated with the control $u_2(\cdot)$ and $\gamma_1$ enjoys the same properties of $\gamma$.
\end{itemize}
We remark (see \eqref{7lboundC}) that $\gamma_1(\tau)$ is
bounded away from $0$ as $\tau\rightarrow 0^{+}$. Moreover, one can see that $z(t)=y_2(t)-y_1(t)$. Therefore
\begin{eqnarray}
\langle\zeta,x_2-x_1 \rangle &=&\langle\zeta, x(\tau)-z(\tau)\rangle + \langle\zeta, z(\tau)\rangle\nonumber\\
&\leq & \|x(\tau)-z(\tau)\|+ \langle\zeta, z(\tau)\rangle.\label{zetax2x1}
\end{eqnarray}
Recalling  (\ref{75dist3}) and (\ref{7linaux}), we obtain 
\begin{equation}\label{7Strictly-convex}
\langle\zeta,x_2-x_1 \rangle\leq (L_4\tau-\gamma_1(\tau)) \Big(\int_{0}^{\tau}|w(s)|ds\Big)^2.
\end{equation}
Thus if $\tau\leq \frac{1}{2L_4}\liminf_{\tau\to 0^+}\gamma (\tau)=:\mathcal{T}_1$ then 
\[
\langle\zeta,x_2-x_1 \rangle\leq -\frac{\gamma_1(\tau)}{2} \Big(\int_{0}^{\tau}|w(s)|ds\Big)^2.
\]
From this inequality the uniqueness of the control steering the origin to $x_1$ in time $\tau$ follows immediately by contradiction.
Setting $\mathcal{T}=\min\lbrace{\mathcal{T}_0,\mathcal{T}_1\rbrace}$ and recalling (\ref{75dist2}) we obtain
(\ref{7final-est}). The proof of the strict convexity is completed by applying Proposition \ref{7convex}.\\
\indent The proof of \eqref{eq:posreach} is entirely analogous, where it suffices to take $\mathcal{T}'=\mathcal{T}_0$.\\
\indent We consider now the statement concerning the ball contained in the reachable set. To this aim, take $u_2\equiv 0$ and set $y_1$
to be the solution of \eqref{7system5} with $u_1$ in place of $w$. Then \eqref{75dist3} yields, for all $t>0$,
\[
 \| x_1(t) - y_1(t)\| \le L_4 t^3.
\]
Recalling \eqref{ballnonaut}, we obtain from the previous inequality that 
\[
\| x_1(t)\| \ge \| y_1(t)\| - L_4t^3\ge \tilde{\gamma} t^2 - L_4 t^3, 
\]
for a suitable constant $\tilde{\gamma}$, which yields in particular that $0$ belongs to the interior of
$\mathcal{R}^t$, $0\le t\le\mathcal{T}$. Since the above argument can be repeated for every point in the
boundary of $\mathcal{R}^t$
and the constant $L_4$ is independent of the reference point, the statement follows by recalling that we
already proved that $\mathcal{R}^t$ is convex for all $0\le t\le\mathcal{T}$.\\
\indent The continuity of $T$ follows easily from the the fact that reachable sets contain a ball (see, e.g.,
Propositions IV.1.2 and IV.1.6 in \cite{BCD}).\\
\indent In the case $M=2$, it suffices to apply the above arguments to each control.
\qed

\medskip

\noindent The following Remark follows immediately from the proof of Theorem \ref{7NL}. 
\begin{Rem}\label{7adjnorm}
Let $x_1(\cdot)$ and $\lambda(\cdot)$ be as in the proof of Theorem \ref{7NL}. For all $0<t\leq\tau=T(x_1)$,
one has $\lambda(t)\in N_{\mathcal{R}^t}(x_1(T(x_1-t)))$. More precisely,
\[
\langle\lambda(t),y-x_1(T(x_1)-t)\rangle\leq -\gamma \|\lambda(t)\|\|y-x_1(T(x_1)-t)\|^2
\]
for all $y\in\mathcal{R}^t$, where $\gamma$ is the constant appearing in \eqref{7final-est}.
\end{Rem}
In fact, put $\lambda (t)$ in place of $\zeta$ in \eqref{zetax2x1}. Then the first summand can be estimated in the same way,
while the upper bound on the second summand, namely the analogue of \eqref{7linaux}, can be obtained through
the same arguments leading to \eqref{7linaux}.      \qed

\medskip

\noindent We conclude this section with a counterexample showing the sharpness of assumption (iii) in Theorem \ref{7NL}.
\begin{Rem}\label{exnonstrict}
An example of a two dimensional nonlinear control system satisfying assumptions (i) and (ii) of Theorem \ref{7NL}
such that the reachable set $\mathcal{R}^\tau$ is not convex for all $\tau > 0$.
\end{Rem}
Consider the control system
\begin{equation}\label{sysexample}
 \begin{cases}
  \dot{x}_1 &= x_2 (1+u)\\
  \dot{x}_2 & = u,\qquad u\in [-1,1].
 \end{cases}
\end{equation}
Set $F(x_1,x_2)=(x_2,0)$, $G(x_1,x_2)=(x_2,1)$, and observe that the assumptions (i) and (ii) of Theorem \ref{7NL} are satisfied, while (iii) is not.
The Hamiltonian for this system is
\[
 \mathcal{H}\big( (x_1,x_2),(p_1,p_2),u \big)=x_2p_1+ u (p_1x_2+p_2),
\]
and Pontryagin's Maximum Principle states that if $\bar{x}(\cdot)$ is an optimal trajectory corresponding to the control $\bar{u}(\cdot)$, then there exists
a function $\lambda=(\lambda_1,\lambda_2)$, never vanishing, and a constant $\lambda_0\le 0$ such that, for a.e. $t$,
\begin{eqnarray*}
 \lambda_1(t)&\equiv&\lambda_1\neq 0 \\
\dot{\lambda}_2 (t)&=&-(1+\bar{u}(t))\lambda_1\\
\bar{x}_2(t)\lambda_1+\bar{u}(t)(\lambda_1\bar{x}_2(t)+\lambda_2(t))+\lambda_0&=&0\\
\bar{u}(t)(\lambda_1\bar{x}_2(t)+\lambda_2(t))&=&\max_{|u|\le 1} u (\lambda_1 \bar{x}_2 (t) + \lambda_2(t)).
\end{eqnarray*}
Now, observe that 
\begin{eqnarray*}
\lambda_1\bar{x}_2(t) + \lambda_2(t)&=&\lambda_1\int_0^t\bar{u}(s)\, ds-\lambda_1 t - \lambda_1\int_0^t\bar{u}(s)\, ds
         +\lambda_2(0)\\
&=&-\lambda_1 t + \lambda_2(0).
\end{eqnarray*}
Therefore the function $\lambda_1\bar{x}_2(t) + \lambda_2(t)$ has at most one zero, so that the optimal control $\bar{u}$ is unique, bang-bang,
and has at most one switching.\\
Fix now $\tau >0$ and $0<s<1$ and consider the control
\[
u_s(t)= 
\begin{cases}
  1& 0<t<s\tau\\
-1 &s\tau < t <\tau.
 \end{cases}
\]
The trajectory of \eqref{sysexample} emanating from the origin and corresponding to the control $u_2$ is, at time $\tau$,
\[
 \bar{x}_s^1 := x_s^1(\tau) = s^2\tau^2,\quad \bar{x}_s^2 := x_s^2(\tau) = \tau (2s-1).
\]
Simple computations show that any other bang-bang control with at most one switching cannot reach $(\bar{x}^1_s,\bar{x}_s^2)$ at a time
$\tau' <\tau$, for all $0<s<1$. Thus $u_s$ is optimal. In particular,
the curve $\gamma(s):=(\bar{x}^1_s,\bar{x}^2_s)$ belongs to the boundary of the reachable set $\mathcal{R}^\tau$. The unique unit normal to the curve $\gamma(s)$
at $s=1/2$ which points outside $\mathcal{R}^\tau$ is 
\[
\zeta = \frac{(2,-\tau)}{\sqrt{4+\tau^2}}.
\]
We compute:
\[
 \big\langle \zeta , \gamma(s)-\gamma(\frac12)\big\rangle = \frac{2\tau^2 (s-\frac12)^2}{\sqrt{4+\tau^2}},
\]
which implies that $\mathcal{R}^\tau$ is not convex. Observe however that $\mathcal{R}^\tau$ has positive reach at $\gamma(\frac12)$ for every $\tau >0$, since
\[
  \big\langle \zeta , \gamma(s)-\gamma(\frac12)\big\rangle \le \frac{8}{(16 + \tau^4)\sqrt{4+\tau^2}}
  \Big\| \gamma(s)-\gamma(\frac12)\Big\|^2.
\]
\qed
\begin{Rem}\label{dimrestr2}
On the assumption $N=2$. 
\end{Rem}
Motivations for the restriction $N=2$ are twofold. First of all, our analysis is based on the switching function of the nonautonomous
system obtained by linearizing around an optimal trajectory, and this method requires $N\le 3$ (see Remark \ref{dimrestr1}). 
Second, the distance between trajectories of the nonlinear system \eqref{7system4} and of the linearized system \eqref{7system5}
is of order two with respect to the control (see \eqref{75dist3}), and this quadratic perturbation can be balanced by
the strict convexity of the reachable set of the linearized system only if $N=2$ (see \eqref{7Strictly-convex}).\qed
\section{Further results for the nonlinear two dimensional case}\label{sec:further}
This section is devoted to proving that the epigraph of the minimum time function has positive reach, under the assumptions of Theorem \ref{7NL}.
To this aim, a results of optimal points, i.e., on points which are \textit{crossed} by an optimal trajectory, is needed.
\subsection{Optimal points}\quad\\
The classical definition of optimal point reads as follows.
\begin{Definition}
Let $x\in\mathbb{R}^N\backslash\lbrace{0\rbrace}$. We say that $x$ is optimal if and only if there exists a point $x_1$
such that $T(x_1)>T(x)$ and a control $u$ with the property that $y^{x_1,u}(\cdot)$ steers $x_1$ to $0$ in the
optimal time $T(x_1)$ and $x=y^{x_1,u}(T(x_1)-T(x))$.
\end{Definition}
The following is the result on optimal points which will be used in the next subsection in order to ensure the positive reach of the
epigraph of the minimum time function. It is based on the same
estimates which lead to the strict convexity of the reachable set, and so it is restricted to two dimensional control systems.
\begin{Theorem}\label{7optimal} Let $N=2$ and let the assumptions of Theorem \ref{7NL} be satisfied. Let $\mathcal{T}>0$
be such that, according to Theorem \ref{7NL}, for all $0\leq\tau<\mathcal{T}$, the reachable set $\mathcal{R}^{\tau}$ 
of \eqref{7system4} satisfies \eqref{7final-est} for all $0<\tau <\mathcal{T}$.
Let $\bar{x}$ be such that $T(\bar{x})<\mathcal{T}$. Then $\bar{x}$ is an optimal point.
\end{Theorem}
\textbf{Proof.} We consider first the case where $G$ is a vector and the control $u$ is one-dimensional. Set $\tau=T(\bar{x})$ and
let $\bar{u}(\cdot)$ be the admissible control steering $\bar{x}$ to $0$ in the optimal time $\tau$, together with the associate
trajectory $\bar{x}(\cdot)$. Set, for all $t\in [0,\tau]$,
\[
A(t)=DF(\bar{x}(t))+DG(\bar{x}(t))\bar{u}(t),\qquad b(t)=G(\bar{x}(t))
\]  
We assume preliminarily that $\bar{x}$ belongs to the boundary of $\mathcal{R}^\tau$ and let, by the Maximum Principle, $\lambda$ be a solution of 
\begin{equation}
\left\{\begin{array}{ll}
\dot{\lambda}(t)=-\lambda(t)A(t),\\
\lambda(\tau)=\zeta,
\end{array}\right.
\end{equation}
where $\zeta\in N_{\mathcal{R}^{\tau}}(\bar{x})$, $\zeta\neq 0$, and, for a.e. $t\in [0,\tau]$,
\begin{equation}\label{usign}
\bar{u}(t)=\mathrm{sign}\langle\lambda(t),b(t)\rangle.
\end{equation}
Set, for $t\in [0,\tau]$,
\[
g(t)=\langle\lambda(t),b(t)\rangle.
\]
We are now going to extend $\bar{u}(\cdot)$ in an interval $[\tau,\tau+\delta]$ for a suitable $\delta>0$, with the property that
the extended control and its associate trajectory satisfy the Maximum Principle.\\
\indent Three cases may occur:
\begin{enumerate}
\item [(i)] $g(\tau)>0$,
\item [(ii)] $g(\tau)<0$,
\item [(iii)] $g(\tau)=0$.
\end{enumerate}
In the first case, we set $\bar{u}(t)=1$ for all $t>\tau$ and let $\bar{x}(\cdot)$ be
the associate trajectory satisfying $\bar{x}(\tau)=\bar{x}$.
We extend analogously $A(\cdot)$, $b(\cdot)$, $\lambda(\cdot)$ and $g(\cdot)$ for $t>\tau$. Set $g(\tau):=\delta_1$. 
Observe that $g$ is locally Lipschitz,
so that, for $t>\tau$,
\[
g(t)=g(\tau)+g(t)-g(\tau)>\delta_1-L_1(t-\tau)
\]
for a suitable constant $L_1$. Therefore we can find $\delta>0$ such that $0\leq\tau+\delta<\mathcal{T}$ and $g(t)>0$ for all $t\in [\tau,\tau+\delta]$, i.e.,
\[
\bar{u}(t)=\mathrm{sign}\ g(t)\quad\forall t\in [\tau,\tau+\delta].
\]
The second case is entirely analogous, by substituting $1$ with $-1$.\\
We consider now the third case. Let the $I_0$, $I_1$ be given by Lemma \ref{7gg0} for the function $g$ in the interval $[0,\tau]$.
Observe that necessarily $\tau\in I_1$, so that, in particular, $g'(\tau)\neq 0$. We set, for $t>\tau$
\[
\bar{u}(t)=1\quad\mathrm{if}\quad g'(\tau)>0
\]
or 
\[
\bar{u}(t)=-1\quad\mathrm{if}\quad g'(\tau)<0
\]
and let $\bar{x}(\cdot)$ be the associate trajectory satisfying $\bar{x}(\tau)=\bar{x}$. We extend analogously
$A(\cdot),b(\cdot),\lambda(\cdot)$ and $g(\cdot)$ for $t>\tau$. Observe that 
\[
\dot{g}(t)=\langle\lambda(t),[F,G](\bar{x}(t))\rangle,
\]
where $[F,G](x)=DG(x)F(x)-DF(x)G(x)$ denotes the Lie bracket of $F$ and $G$.
Therefore $\dot{g}$ is continuous, so that there exists $\delta>0$ such that the sign of $\dot{g}(t)$ equals the sign of
$\dot{g}(\tau)$ for all $t\in [\tau-\delta,\tau+\delta]$. Therefore our construction of $\bar{u}(\cdot)$ on
$[0,\tau+\delta]$ is such that for a.e. $t\in [0,\tau+\delta]$,
\[
\bar{u}(t)=\text{sign}\ g(t).
\] 
Consequently, all conclusions of Theorem \ref{7NL} hold up to the time $\tau+\delta$. In particular,
for all $t\in [0,\tau+\delta]$, $\bar{x}(t)\in\mathrm{bdry}\,\mathcal{R}^{t}$. Since $T(\cdot)$ is continuous in a neighborhood of the trajectory
$\bar{x}(\cdot)$, we obtain that $\bar{u}(\cdot)$ steers the origin optimally to $\bar{x}(\tau+\delta)$ in time $\tau+\delta$.
Since the above argument can be applied also to the reversed dynamics $\dot{x}=-F(x)+G(x)u$, $u\in [-1,1]$, then $T(\bar{x}(\tau+\delta))=\tau+\delta$.\\
Let us now drop the assumption $\bar{x}\in\mathrm{bdry}\,\mathcal{R}^\tau$. Since $T$ is strictly decreasing along the optimal trajectory $\bar{x}(\cdot)$,
and so $\bar{x}(t)\in\mathrm{bdry}\,\mathcal{R}^{\tau-t}$ for all $0<t<\tau$, there exists a nontrivial adjoint vector $\lambda(\cdot)$ which uniquely determines
$\bar{u}(t)$ as in \eqref{usign} up to the time $\tau$. Thus the above argument can be applied also to $\bar{x}$.\\
If $G$ is a $2\times 2$ matrix, it suffices to perform the above construction for each column of $G$. The proof is concluded.\qed

\medskip

We conclude the section with two corollaries. The first one is an immediate consequence of the proof of Theorem \ref{7optimal}.
\begin{Corollary}\label{7extoptimal}
Under the same assumptions of Theorem \ref{7optimal}, let $\tau=T(\bar{x})<\tau_1<\mathcal{T}$. Then there exists $x_1\in\mathrm{bdry}\,\mathcal{R}^{\tau_1}$
and a control $u_1:[\tau,\tau_1]\rightarrow [-1,1]$ such that the trajectory $\tilde{x}(\cdot)$ corresponding to the control
\begin{equation*}
\tilde{u}(t)=\left\{\begin{array}{ll}
\bar{u}(t)\quad 0\leq t\leq \tau,\\
u_1(t)\quad \tau<t\leq\tau_1.
\end{array}\right.
\end{equation*}
and such that $\tilde{x}(0)=x_1$ reaches $0$ in the optimal time $\tau_1$ and moreover $\tilde{x}(\tau_1-\tau)=\bar{x}$.
\end{Corollary}
From Corollary \ref{7extoptimal} we obtain that $\bar{x}(\cdot)$ and $\lambda(\cdot)$ in the proof of Theorem \ref{7optimal} can
be extended up to the time $\mathcal{T}$. Therefore, we obtain the following further corollary.
\begin{Corollary}\label{7Maxhal}
Under the assumptions of Theorem \ref{7NL}, the maximized Hamiltonian along $\bar{x}(\cdot)$ associated
with $\lambda(\cdot)$ is constant in $[0,\mathcal{T})$, i.e.,
 \[
 H(\bar{x}(t),\lambda(t))=C\quad\forall t\in [0,\mathcal{T}). 
 \] 
\end{Corollary}
\textbf{Proof.} Let $G$ be a vector and so the control $u$ be scalar. Then
\[
H(\bar{x}(t),\lambda(t))=\langle\lambda(t),F(\bar{x}(t))\rangle+|\langle\lambda(t),G(\bar{x}(t))\rangle|.
\]
Observe that the switching function $g(t)=\langle \lambda (t),G(\bar{x}(t))\rangle$ vanishes at most on a countable subset of $[0,\mathcal{T})$.
Therefore, for a.e. $t\in [0,\mathcal{T})$, we have
\[
\frac{d}{dt}H(\bar{x}(t),\lambda(t))=\Big(-\langle\lambda(t),[F,G](\bar{x}(t))\rangle+\langle\lambda(t),[F,G](\bar{x}(t))\rangle\Big)
\mathrm{sign}\, g(t)=0.
\]
If $G$ is a $2\times 2$ matrix, it suffices to perform the above computation for each column. The proof is concluded. \qed
\subsection{The epigraph of the minimum time function has positive reach}\quad\\
The present section is devoted to studying the ``convexity type'' of the minimum time function $T(\cdot)$, in the case where the dynamics satisfies
a weak controllability condition, i.e., the function $T(\cdot)$ is merely continuous. The statement is two dimensional, since it is based on
Theorem \ref{7optimal}.
\begin{Theorem}\label{7posreach}
Let $N=2$ and let the assumptions of Theorem \ref{7NL} hold. Let $\mathcal{T}$ be given by Theorem \ref{7NL}.
Then for every $0<\tau<\mathcal{T}$ the epigraph of the minimum time function $T(\cdot)$ on $\mathcal{R}^{\tau}$ has positive reach.
\end{Theorem}
\begin{Corollary}
Under the same assumptions of Theorem \ref{7posreach} the minimum time function $T$ satisfies all the properties listed in Theorem \ref{diff}.
\end{Corollary}
Before beginning the proof of Theorem \ref{7posreach} we introduce the minimized Hamiltonian and study its sign.
\begin{Definition}
Let $x,\zeta\in\mathbb{R}^N$. We define the minimized Hamiltonian for the control system in \eqref{7system4} as 
\[
h(x,\zeta)=\langle\zeta,F(x)\rangle+\min_{u\in\mathcal{U}}\langle\zeta,G(x)u\rangle.
\]
\end{Definition}
\begin{Proposition}\label{7signh}
Let $x$ belong to the boundary of the reachable set $\mathcal{R}^{\tau}$
for \eqref{7system4} for some $\tau>0$. Let $\zeta\in N^F_{\mathcal{R}^\tau}(x)$\footnote{here $N^F_{\mathcal{R}^\tau}(x)$ denotes the Fr\'echet
normal cone to $\mathcal{R}^\tau$ at $x$, i.e., all vectors $v$ such that $\limsup_{\mathcal{R}^\tau\ni y\to x}\langle v , (y-x)/\| y-x\|\rangle\le 0$}.
Then $h(x,\zeta)\leq 0$.
\end{Proposition}
\textbf{Proof.} Let $\bar{u}(\cdot)$ be an admissible control steering $x$ to $0$ in time $\tau$, together with
the associate trajectory $\bar{x}(\cdot)$. Then, for all $0\leq t\leq\tau$ the point $\bar{x}(t)$ belongs to $\mathcal{R}^{\tau}$, so
that, by definition of Fr\'echet normal we have
\[
\limsup_{t\rightarrow 0^+}\ \Big\langle\zeta,\frac{\bar{x}(t)-x}{\| \bar{x}(t)-x\|}\Big\rangle\leq 0.
\]
\indent Observing that $\|x(t)-x\|\leq Kt$ for a suitable constant $K$, we have
\[
\limsup_{t\rightarrow 0}\ \Big\langle\zeta,\frac{\bar{x}(t)-x}{t}\Big\rangle\leq 0.
\]
In other words,
\begin{eqnarray*}
0&\geq &\limsup_{t\rightarrow 0}\Big\langle\zeta,\frac{1}{t}\int_{0}^t\big(F(\bar{x}(s))+G(\bar{x}(s))\bar{u}(s)\big)ds\Big\rangle\\
&=&\langle\zeta ,F(x)\rangle +\limsup_{t\rightarrow 0}\Big\langle\zeta, G(x)\frac{\int_0^t\bar{u}(s)ds}{t}\Big\rangle.
\end{eqnarray*}
Let $t_n\rightarrow 0$ be a sequence such that $\lim_{n\rightarrow\infty}\frac{1}{t_n}\int_0^{t_n}\bar{u}(s)ds:=\tilde{u}$ exists.
By the convexity of $\mathcal{U}$, $\tilde{u}\in\mathcal{U}$, and so $h(x,\zeta)\le \langle\zeta,F(x)\rangle+\langle\zeta,G(x)\tilde{u}\rangle\le 0$.
\qed

\medskip

We are now ready to prove Theorem \ref{7posreach}.

\medskip

\noindent\textbf{Proof of Theorem \ref{7posreach}.} Let $x\neq 0$ be such that $T(x)< \mathcal{T}$ and let $(\bar{u}(\cdot),\bar{x}(\cdot))$
be an optimal pair for $x$. By Maximum Principle\footnote{observe that we are applying Theorem \ref{th:PMPnonlin} for the \textit{reversed dynamics} $\dot{x}=-F(x)-G(x)u$} there exists $0\neq\zeta\in N_{\mathcal{R}^{T(x)}}(x)$
such that the adjoint arc $\lambda$, with 
\begin{equation*}
\left\{\begin{array}{ll}
\dot{\lambda}(t)=\lambda(t)\big(DF(\bar{x}(t))+DG(\bar{x}(t))\bar{u}(t)\big),\\
\lambda(T(x))=\zeta
\end{array}\right.
\end{equation*}
satisfies
\[
\big\langle\lambda(t),F(\bar{x}(t))+G(\bar{x}(t))\bar{u}(t)\big\rangle=h(\bar{x}(t),\lambda)\quad\text{for a.e. }t\in (0,T(x)).
\]
We claim that 
\begin{equation}\label{7normepi}
(\zeta,h(x,\zeta))\in N^P_{\mathrm{epi}(T)}(x,T(x)),
\end{equation}
i.e., there exists a constant $\sigma>0$ such that, for all $y\in\mathbb{R}^N$ with $0<T(y)<\mathcal{T}$ and for all $\beta\geq T(y)$,
we have
\begin{equation}\label{7normepi2}
\big\langle(\zeta,\theta),(y,\beta)-(x,T(x))\big\rangle\leq\sigma \|(\zeta,\theta)\|\Big(\|y-x\|^2+|\beta -T(x)|\Big),
\end{equation}
where $\theta=h(x,\zeta)$, and, moreover,
\begin{eqnarray}\label{7indep}
\sigma\ \text{is independent of}\ x\ \text{and}\ \zeta.
\end{eqnarray}
Indeed, we consider two cases:
\begin{enumerate}
\item [(a)] $T(y)\leq T(x)$;
\item [(b)] $T(y)>T(x)$. 
\end{enumerate}
In the first case, $y\in\mathcal{R}^{T(x)}$, so that by Theorem \ref{7NL}
\[
\langle\zeta,y-x\rangle\leq 0. 
\]
If $\beta\geq T(x)$ then (\ref{7normepi2}) is automatically satisfied, since $\theta\leq 0$ by Proposition \ref{7signh}.
If instead $\beta<T(x)$, we set $x_1=\bar{x}(\beta)$.\\
 \indent We estimate first $\langle\zeta,y-x_1\rangle$. Since $y\in\mathcal{R}^{\beta}$, recalling Remark \ref{7adjnorm},
we have for suitable constants $K_1,K_2$ given by Gronwall's Lemma, 
 \begin{eqnarray*}
 \langle\zeta,y-x_1\rangle &=& \langle\lambda(\beta),y-x_1\rangle +\langle\lambda(T(x))-\lambda(\beta),y-x_1\rangle\\
 &\leq & \langle\lambda(T(x))-\lambda(\beta),y-x_1\rangle \leq  K_1\|\lambda(T(x))\|\, |T(x)-\beta|\, \|y-x_1\|\\
 &\leq &  K_1\|\lambda(T(x))\|\, |T(x)-\beta| \big(\|y-x\|+\|x_1-x\|\big)\\
 &\leq & K_1\|\zeta\|\, |T(x)-\beta|\,\big(\|y-x\|+K_2|T(x)-\beta|\big)\\
 &\leq & K_3 \|\zeta\|\big(\|y-x\|^2+|T(x)-\beta|^2\big)
 \end{eqnarray*}
 for another suitable constant $K_3$.\\
 \indent Second, we estimate $\langle\zeta,x_1-x\rangle$. We have 
 \begin{eqnarray*}
 \langle\zeta,x_1-x\rangle &=&\int_\beta^{T(x)}\big\langle\lambda(T(x)),F(\bar{x}(s))+G(\bar{x}(s))\bar{u}(s)\big\rangle\, ds\\
 &=&\int_\beta^{T(x)}\big\langle\lambda(s),F(\bar{x}(s))+G(\bar{x}(s))\bar{u}(s)\big\rangle\, ds\\
 &+&\int_\beta^{T(x)}\big\langle\lambda(T(x))-\lambda(s),F(\bar{x}(s))+G(\bar{x}(x))\bar{u}(s)\big\rangle\, ds\\
 &\leq &(T(x)-\beta)h(x,\zeta)+K_4\|\zeta\||T(x)-\beta|^2,
 \end{eqnarray*}
 for a suitable constant $K_4$, owing to Corollary \ref{7Maxhal} (applied to the reversed dynamics $\dot{x}=-F(x)-G(x)u$).
 Therefore, since $h(x\zeta)\le 0$,
 \[
 \langle(\zeta,\theta),(y,\beta)-(x,T(x))\rangle\leq (K_3+K_4)\|\zeta\|\big(\|y-x\|^2+|T(x)-\beta|^2\big),
  \]
and the proof for the case (a) is concluded by observing that $K_3$ and $K_4$ are independent of $\zeta$ and $x$.

In the second case we need to use the optimality of $x$. We observe first that, since $\theta\leq 0$, we only need
to prove (\ref{7normepi2}) for $\beta=T(y)$. Recalling Corollary \ref{7extoptimal}, we can extend the control $\bar{u}$ up to
the time $T(y)$ so that the associated trajectory (still denoted by $\bar{x}(\cdot)$) remains optimal. Let also $\lambda$ be the
extended adjoint vector and denote by $\tilde{x}(\cdot)$ the trajectory of the reversed dynamics associated with the extended control $\bar{u}$, i.e.,
 \begin{equation*}
\left\{\begin{array}{ll}
\dot{\tilde{x}}(t)=-F(\tilde{x}(t))-G(\tilde{x}(t))\tilde{u}(t),\\
\tilde{x}(0)=0,
\end{array}\right.
\end{equation*}    
where $\tilde{u}(t)=\bar{u}(T(y)-t)$.\\
Set $x_1=\tilde{x}(T(y))$. We estimate first $\langle\zeta,y-x_1\rangle$. We have, by arguing similarly as before,
\begin{eqnarray*}
\langle\zeta,y-x_1\rangle &=&\langle\lambda(T(y)),y-x_1\rangle +\langle\lambda(T(x))-\lambda(T(y)),y-x_1\rangle\\
& & (\text{the first summand is $\le 0$ by the construction in Theorem \ref{7optimal}})\\
&\leq & \langle\lambda(T(x))-\lambda(T(y)),y-x_1\rangle\\
&\leq & K_5 \|\zeta\|\big(|T(y)-T(x)|^2+\|y-x_1\|^2\big).
\end{eqnarray*}
On the other hand,
\begin{eqnarray*}
\langle\zeta,x_1-x\rangle &=&\int_{T(x)}^{T(y)}\big\langle\zeta,-F(\tilde{x}(s))-G(\tilde{x}(s))\tilde{u}(s)\big\rangle ds\\
&=&\int_{T(x)}^{T(y)}\big\langle\lambda(s),-F(\tilde{x}(s))-G(\tilde{x}(s))\tilde{u}(s)\big\rangle ds\\
&+&\int_{T(x)}^{T(y)}\big\langle\lambda(T(x))-\lambda(s),-F(\tilde{x}(s))-G(\tilde{x}(s))\tilde{u}(s)\big\rangle ds\\
&\leq &\int_{T(x)}^{T(y)}\max_{u\in\mathcal{U}}\big\langle\lambda(s),-F(\tilde{x}(s))-G(\tilde{x}(s))u\big\rangle ds+K_6\|\zeta\|(T(y)-T(x))^2\\
& & (\text{for a suitable constant}\ K_6\ \text{given by Gronwall's Lemma}).
\end{eqnarray*}

Recalling Corollary \ref{7Maxhal}, the maximized Hamiltonian in the integral of the first summand is constant. Therefore we obtain
\[
\langle\zeta,x_1-x\rangle\leq -h(x,\zeta)(T(y)-T(x))+K_6\|\zeta\|\, |T(y)-T(x)|^2.
\]
Combining the above estimates we obtain finally
\[
\langle(\zeta,\theta),(y,T(y))-(x,T(x))\rangle\leq (K_5+K_6)\|\zeta\|\big(\|y-x\|^2+|T(y)-T(x)|^2\big),
\]
and the proof of the claim is concluded, by observing, again, that $K_5,K_6$ are independent of $x$ and $\zeta$.\\
\indent In order to conclude the proof we observe that $N^P_{\mathrm{epi}(T)}$ is pointed at every point $(x,T(x))$, $x\in\mathcal{R}^{\tau}$,
since it is easy to see that the projection of every $(\zeta,\theta)\in N^P_{\mathrm{epi}(T)}(x,T(x))$ onto $\mathbb{R}^N$ is normal to the strictly convex set $\mathcal{R}^\tau$.
Therefore, we can apply Corollary 3.1 in \cite{K}, with $\Omega_P=\mathrm{int}\mathcal{R}^{\tau}$, which shows that $\mathrm{epi}(T)$ has positive reach.
\qed
\section{Appendix}\label{sec:app}
This section is devoted to some technical lemmas which are used in the proof of the main results.
\begin{Lemma}\label{7app-intineq} 
Let $K:(a,b)\rightarrow [0,1]$ be measurable and let $k\in\mathbb{N}$. Then 
\[
\int_a^b(t-a)^kK(t)dt\geq\frac{1}{k+1}\Big(\int_a^bK(t)dt\Big)^{k+1},
\]
and 
\[
\int_{a}^{b}(b-t)^kK(t)dt\geq\frac{1}{k+1}\Big(\int_a^bK(t)dt\Big)^{k+1}.
\]
\end{Lemma}
\textbf{Proof.}
Indeed,
\[
\int_{a}^{b}(t-a)^kK(t)dt=k!\int_a^b\int_{t_{k}}^b\ldots \int_{t_1}^bK(t_0)dt_0\ldots dt_k.
\]
Since $K(t)\in [0,1]$ for a.e. $t\in [0,1]$, we obtain that 
\[
\int_{a}^{b}(t-a)^kK(t)dt\ge k!\int_a^bK(t_k)\int_{t_{k}}^bK(t_{k-1})\ldots \int_{t_1}^bK(t_0)dt_0\ldots dt_k
\] 
One can easily prove that
\[
\int_a^bK(t_k)\int_{t_{k}}^bK(t_{k-1})\ldots \int_{t_1}^bK(t_0)dt_0\ldots dt_k=\frac{1}{(k+1)!}\Big(\int_a^bK(t)dt\Big)^{k+1}
\]
The proof is complete. \qed
\begin{Lemma}\label{7ineqg}
Let $f:[a,b] \rightarrow \mathbb{R}$ be absolutely continuous and monotone. Fix $k\geq 0$ and assume that there exists
$C\in\mathbb{R}^+$ such that 
\begin{equation}\label{7bound a}
|f'(s)|\geq C(s-a)^k\quad\forall s\in [a,b].
\end{equation}
Then, either $f$ has no zeros in $(a,b)$ and then, for all $s\in (a,b)$ either
\[
|f(s)|\geq\frac{C}{k+1}(s-a)^{k+1}\ \ \text{or}\ \ |f(s)|\geq\frac{C}{k+1}(b-s)^{k+1},
\]
or there exists $c\in (a,b)$ such that $f(c)=0$ and then, for all $s\in [a,b]$,
\[
|f(s)|\geq \frac{C}{k+1}|c-s|^{k+1}
\]
The same conclusions hold if (\ref{7bound a}) is substituted by 
\[
|f'(s)|\geq C(b-s)^{k}\quad\forall s\in [a,b].
\]
\end{Lemma}
\textbf{Proof.}
We treat the case where $f$ is nondecreasing, while the other
one can be handled by taking $-f$. So (\ref{7bound a}) now reads as
\[
f'(s)\geq C(s-a)^k\quad\forall s\in [a,b].
\]
If $f$ has no zeros, we have two cases, namely $f(s)>0$ for all $s\in (a,b)$ or $f(s)<0$ for all $s\in (a,b)$. For the first case 
\[
f(s)-f(a)=\int_{a}^sf'(t)dt\geq C\int_a^s(t-a)^kdt=\frac{C}{k+1}(s-a)^{k+1},
\]
so that
\[
f(s)\geq \frac{C}{k+1}(s-a)^{k+1}.
\]
In the second case, 
\begin{eqnarray*}
f(b)-f(s)&=&\int_s^bf'(t)dt\geq C\int_{s}^b(t-a)^kdt\\
&=&\frac{C}{k+1}\big[(b-a)^{k+1}-(s-a)^{k+1}\big]\geq\frac{C}{k+1}(b-s)^{k+1}.
\end{eqnarray*}
Therefore,
\[
f(s)\leq f(b)-\frac{C}{k+1}(b-s)^{k+1}\leq -\frac{C}{k+1}(b-s)^{k+1}.
\]
Assume now that there exists $c\in (a,b)$ such that $f(c)=0$. Then, for all $s\in [a,c]$ we have 
\[
-f(s)=f(c)-f(s)=\int_s^cf'(t)dt\geq\frac{C}{k+1}(c-s)^{k+1},
\]
while for all $s\in [c,b]$ we have 
\begin{eqnarray*}
f(s)&=&f(s)-f(c)=\int_c^sf'(t)dt\geq\frac{C}{k+1}\big[(s-a)^{k+1}-(c-a)^{k+1}\big] \\
&\geq &\frac{C}{k+1}(s-c)^{k+1},
\end{eqnarray*}
and the proof is complete. \qed
\begin{Proposition}\label{7convex}
Let $K\subset\mathbb{R}^N$ be compact and assume that there exist $\gamma>0$ and $p>1$ with the following property:
for every $x\in\mathrm{bdry}\, K$, there exists $\zeta\neq 0$ such that for every $y\in K$ one has
\begin{equation}\label{7strict}
\langle\zeta,y-x\rangle\leq -\gamma\|\zeta\|\|y-x\|^p.
\end{equation}
Then $K$ is convex (with nonempty interior) and, for each $x\in\mathrm{bdry}\, K$, (\ref{7strict}) is satisfied by all $\zeta\in N_K(x)$.
\end{Proposition}
\textbf{Proof.} We show first that $K$ is strictly convex. To this aim, assume by contradiction that there exist points $x_1\neq x_2\in K$
such that the segment $]x_1,x_2[$ is not contained in the interior of $K$. Let $0<t<1$ be such that $x_t=(1-t)x_1+tx_2\in\mathrm{bdry}\, K$ and let $\zeta\neq 0$
be such that (\ref{7strict}) holds with $x_t$ in place of $x$. Obviously, $\langle\zeta, x_2-x_1\rangle=0$, and this is a contradiction.
Since $K$ is convex with nonempty interior, for each $x\in\mathrm{bdry}\, K$
the normal cone $N_K(x)$ is pointed and so it is the convex hull of its exposed rays (see \cite{PoR}).\\
\indent We now use Theorem 4.6 in \cite{CM} and see that for every unit vector $w$ belonging to an exposed ray of $N_K(x)$,
there exists a sequence $x_n\rightarrow x$ such that 
\[
N_K(x_n)=\mathbb{R}^+w_n,\; \|w_n\|=1\quad\text{and } w_n\to w.
\]
Of course (\ref{7strict}) holds with $x_n$ (resp., $w_n$) in place of $x$ (resp., $\zeta$), so that by passing to the limit,
$w$ also satisfies (\ref{7strict}). By taking convex combinations, we conclude the proof.
\qed


\begin{thebibliography}{999999}

\bibitem{BCD} {\sc M. Bardi \& I. Capuzzo-Dolcetta}, {\em Optimal Control and
Viscosity Solutions of Hamilton--Jacobi--Bellman Equations}, Birkh\"auser, Boston (1997).

\bibitem{BP} {\sc U. Boscain \& B. Piccoli},  {\em Optimal syntheses for
control systems on 2-D manifolds}, Springer, Berlin (2004).

\bibitem{cafra} {\sc P. Cannarsa, \& H. Frankowska}, Interior sphere property of attainable sets and time
optimal control problems, ESAIM Control Optim. Cal. Var. 12 (2006), 350-370. 

\bibitem{camapw} {\sc P. Cannarsa, F. Marino \& P. R. Wolenski}, {\em On the minimum time function for differential inclusions},
Discrete Contin. Dyn. Syst. Ser. A, in print.

\bibitem{CaKh} {\sc P. Cannarsa \& Khai T. Nguyen}, {\em Exterior sphere condition and time optimal control
for differential inclusions}, SIAM J. Control Optim., to appear.

\bibitem{CS0}  {\sc P. Cannarsa \& C. Sinestrari}, {\em Convexity properties of the
minimum time function}, Calc. Var. Partial Differential Equations 3 (1995), 273-298.

\bibitem{CS}{\sc P. Cannarsa \& C. Sinestrari}, {\em Semiconcave Functions, Hamilton-Jacobi Equations, and optimal Control}, Birkh\"auser,
Boston (2004).

\bibitem{cesari} {\sc L. Cesari}, {\em Optimization--Theory and Applications}, Springer, New York (1983).

\bibitem{clarke} {\sc F. H. Clarke} {\em Optimization and Nonsmooth Analysis}, SIAM, Philadelphia (1990).

\bibitem{CLSW} {\sc F. H. Clarke, Yu. S. Ledyaev, R. J. Stern \& P. R. Wolenski}, {\em Nonsmooth Analysis and Control Theory},
Springer, New York (1998).

\bibitem{CM} {\sc G. Colombo \& A. Marigonda}, {\em Differentiability properties
for a class of non-convex functions}, Calc. Var. Partial Differential Equations 25 (2006), 1--31.

\bibitem{CMW} {\sc G. Colombo, A. Marigonda \& P. R. Wolenski}, {\em Some new regularity properties for the minimal time function},
SIAM J. Control Optim. 44 (2006), 2285--2299.

\bibitem{CK} {\sc G. Colombo \& Khai T. Nguyen}, {\em On the structure of the Minimum Time Function}, SIAM J. Control Optim. 48 (2010),
4776--4814.

\bibitem{federer} {\sc H. Federer}, {\em Curvature measures}, Trans. Amer. Math. Soc. 93 (1959), 418--491.

\bibitem{hls} {\sc H. Hermes \& J. P. LaSalle}, \emph{Functional analysis and
  time optimal control}, Academic Press, New York-London (1969).

\bibitem{kre} {\sc A. Krener}, {\em A generalization of Chow's theorem and the bang-bang theorem to nonlinear control problems}, SIAM J. Control 12 (1974),
43--52.
  
\bibitem{K} {\sc Khai T. Nguyen}, {\em Hypographs satisfying an external sphere condition and the regularity of the minimum time function},
J. Math. Anal. Appl. 372 (2010), 611--628.

\bibitem{PoR} {\sc R. T. Rockafellar}, {\em Clarke's tangent cones and the boundaries of closed sets in $\mathbb{R}^n$},
Nonlinear Analysis, Theory, Methods and Applications. 3 (1979), 145-154.

\bibitem{suss} {\sc H. Sussmann}, {\em A bang-bang theorem with bounds on the number of switchings}, SIAM J. Control Optim. 17 (1979), 629--651.
\end{thebibliography}
\end{document}